\documentclass{article}
\usepackage{latexsym}
\usepackage{amssymb}
\usepackage{amsmath}
\usepackage{graphicx}
\usepackage{color}
\begin{document}

\newtheorem{thm}{Theorem}[section]
\newtheorem{prop}[thm]{Proposition}
\newtheorem{lem}[thm]{Lemma}
\newtheorem{cor}[thm]{Corollary}
\def\dem{\par\noindent \emph{Proof.} \ \,}
\def\edem{\rm\hfill$\Box$\bigskip}
\def\eq{\operatorname{Eq}\,}
\def\m{^{-1}}
\def\vhi{\varphi}

\newcommand{\cref}[1]{Corollary~\ref{#1}}
\newcommand{\lref}[1]{Lemma~\ref{#1}}
\newcommand{\tref}[1]{Theorem~\ref{#1}}
\newcommand{\secref}[1]{Section~\ref{#1}}

\def\veps{\varepsilon}
\def\ll{{\textstyle \ast}}
\def\rr{{\scriptscriptstyle \triangle}}

\title{Latin bitrades, dissections of equilateral triangles and abelian groups}

\author{\Large 
Ale\v s Dr\'apal\footnote{Supported by grant MSM 0021620839}\\
Department of Mathematics \\ Charles University \\ Sokolovsk\'a 83 \\ 186 75 Praha 8 \\ Czech Republic \\
~ \\
Carlo H\"{a}m\"{a}l\"{a}inen\footnote{Supported by Eduard \v Cech center, grant LC505.}  \\
Department of Mathematics \\ Charles University \\ Sokolovsk\'a 83 \\ 186 75 Praha 8 \\ Czech Republic \\
{\texttt carlo.hamalainen@gmail.com}\\
~ \\
V\'{\i}t\v ezslav Kala\footnote{V. Kala supported by GAUK 8648/2008.} \\
Department of Mathematics \\ Charles University \\ Sokolovsk\'a 83 \\ 186 75 Praha 8 \\ Czech Republic
}

\maketitle

\begin{abstract}
Let $T = (T^{\textstyle \ast}, T^{\scriptscriptstyle \triangle})$ be a spherical latin bitrade. With each 
$a=(a_1,a_2,a_3)\in T^{\textstyle \ast}$ associate a set of linear equations 
$\eq(T,a)$ of the form $b_1+b_2=b_3$, where $b = (b_1,b_2,b_3)$
runs through $T^{\textstyle \ast} \setminus \{a\}$. Assume $a_1 = 0 = a_2$
and $a_3  = 1$. Then $\eq(T,a)$ has in rational numbers a unique solution 
$b_i = \bar b_i$. Suppose that $\bar b_i \ne \bar c_i$ for all 
$b,c \in T^{\textstyle \ast}$
such that $b_i \ne c_i$ and $i \in \{1,2,3\}$. We prove that
then $T^{\scriptscriptstyle \triangle}$ can be interpreted as
a dissection of an equilateral triangle. We also consider
group modifications of latin bitrades and show that the methods
for generating the dissections can be used for a proof that
$T^{\textstyle \ast}$ can be embedded into the operational
table of a finite abelian group, for every spherical
latin bitrade $T$.
\end{abstract}

\section{Introduction}\label{1}

Consider an equilateral triangle $\Sigma$ that is dissected
into a finite number of equilateral triangles. 
Dissections will be always assumed to be nontrivial, and so the
number of dissecting triangles is at least four. 
Denote by $a$, $b$ and $c$
the lines induced by the sides of $\Sigma$. It is easy
to realize that each side of a dissecting triangle has to be parallel
to $a$ or $b$ or $c$. If $X$ is a vertex
of a dissecting triangle, then $X$ is a vertex of exactly one,
three or six dissecting triangles. Suppose that there is
no vertex with six triangles and consider triples 
$(u,v,w)$ of lines that are parallel to $a$, $b$ and $c$, respectively,
and meet in  a vertex of a dissecting triangle that is not a vertex
of $\Sigma$. The set of all these triples together with the triple
$(a,b,c)$ will be denoted by $T^\ll$, and by $T^\rr$ we shall
denote the set of all triples $(u,v,w)$ of lines that are yielded 
by sides of a
dissecting triangle (where $u$, $v$ and $w$ are again parallel to
$a$, $b$ and $c$, respectively). Observe that the following conditions
hold:
\begin{enumerate}
\item[(R1)] Sets $T^\ll$ and $T^\rr$ are disjoint;
\item[(R2)] for all $(p_1,p_2,p_3)\in T^\ll$ and 
all $r,s \in \{1,2,3\}$, $r \ne s$, there exists exactly one
$(q_1,q_2,q_3) \in T^\rr$ with $p_r = q_r$ and $p_s = q_s$; and
\item[(R3)] for all $(q_1,q_2,q_3)\in T^\rr$ and
all $r,s \in \{1,2,3\}$, $r \ne s$, there exists exactly one
$(p_1,p_2,p_3) \in T^\ll$ with $q_r = p_r$ and $q_s = p_s$.
\end{enumerate}

Note that (R2) would not be true if there had existed six dissecting
triangles with a common vertex. Conditions (R1--3) are, in fact, axioms
of a combinatorial object called latin bitrades \cite[p.148]{wanenc}.
This way of their
construction was described in \cite{d1,d2} and here we shall consider
the converse approach, i.~e.~determining when a latin bitrade yields 
a dissection. The topic 
obtained a strong impetus recently when Cavenagh and 
Lison\v ek observed \cite{cl} that spherical latin bitrades
are equivalent to cubic 3-connected bipartite planar graphs \cite{hmk}.
Since there exists a computer package \cite{plantri} that uses
an algorithm of 
Batagelj \cite{bag} for a fast generation of the latter objects, 
it is natural to adapt it for generation of 
dissections. Note that Wanless enumerated small latin bitrades independently
of \cite{plantri}. His results \cite{wan} influenced the
above mentioned discovery of Cavenagh and Lison\v ek.

Some dissections are easy to find, but some are difficult to unravel,
and there seems to exist no efficient way to generate 
dissections directly without resorting to some kind of abstraction.
A spherical latin bitrade has to possess an embedding into a cyclic
group if it can yield a dissection, but it is not clear if every
embedding into a cyclic group can be interpreted as a dissection.

In every dissection of $\Sigma$ there exist exactly three dissecting
triangles with a vertex that is also a vertex of $\Sigma$.
If we remove one or two of these triangles we obtain a pentagon
or a quadrangle. 
Starting from them we can construct further pentagons or
quadrangles by adding a triangle to one of the sides.
Any arising trapezoid quadrangle may be completed to a new
dissection, and we
can regard two dissections
as related when they can be built in this manner from a common origin.
However, the nature of such origins is not well understood yet.
A better understanding might help to explicate the structure
of all spherical latin bitrades and to influence a solution of 
Barnette's conjecture. 

Lines parallel to $a$, $b$ or $c$ that are induced by a side
of a dissecting triangle will be called \emph{dissecting lines}.
For any of them the union of all sides of dissecting triangles
that are incident to the line forms
one or more contiguous segments. If there are two or more segments,
then upon the line there exist 
two dissecting vertices such that all
triangles in between are cut by the line into two parts. If such
a situation arises for no dissecting line and if no
dissecting vertex is incident to six dissecting triangles,
then we call the dissection \emph{separated}. 

Our initial motivation was to find a relatively efficient algorithm
that would produce all separated dissections of small orders.
Experiments indicated that for small orders nearly all spherical latin
bitrades produce a dissection. The question was how to 
systematically reverse the process described in the beginning
of this section. The first step in the reverse process is the choice
of $u = (u_1,u_2,u_3) \in T^\ll$ that induces the sides $a$, $b$
and $c$ of $\Sigma$. With $u$ chosen we convert $T= (T^\ll, T^\rr)$
into a system of $s-1$ linear equations with $s-1$ variables, each
of which corresponds to a dissecting line. The integer $s$ coincides
with the number of dissecting triangles. The system has always
a unique solution and our first main result (\tref{21}) states
that a \emph{pointed} spherical latin bitrade $(T,u)$ determines
a separated dissection if and only if it is true that inside each
group the same value is never assigned to two different unknowns.
By a group we mean here the set of all dissecting lines parallel
to $a$ (or $b$, or $c$).

The proof of \tref{21} is quite long and stretches over Sections~\ref{2}
and~\ref{3}. The result may look to be somewhat counterintuitive
since there seems to be no obvious single reason why a solution 
could not induce a covering of $\varSigma$ by triangles in which at least
one overlapping occurs. We shall use elementary geometry to study
the local structure of the overlapping and show that in fact it 
never occurs.

After proving \tref{21} we realized that some arguments remain valid
in a weaker form even when the solution is not necessarily separated.
These arguments now appear as the three lemmas of \secref{2}.
They indicate how to prove that every spherical latin bitrade can be embedded
into a finite abelian group, our second main result. 
The additional arguments that are needed for the result appear
in \secref{T}.

The connection of dissections 
to counting modulo $n$ becomes clear when we place a dissection
upon a grid formed by equilateral triangles of size $1$ in such a way
that every dissecting triangle becomes a union of basic grid triangles.
If the length of the dissected triangle is $n$, then the dissecting
vertices can be interpreted as cells in the addition table modulo $n$.
This establishes an embedding of $T(\ll)$ into 
$\mathbb Z_n(+)$.

When the solution is not separated, then we do not get an embedding,
but we still get a mapping $T(\ll) \to \mathbb Z_n(+)$ that behaves
like a homomorphism. We shall rather speak about a \emph{homotopy} 
to stress the
fact that the mapping is defined independently for each of the three
groups. \lref{24} states that if $v =(v_1,v_2,v_3)\in T^\ll$ and
$i \in \{1,2,3\}$ are such that $u_i \ne v_i$, then there exists
a homotopy to a finite cyclic group that assigns to $u_i$ and $v_i$
different values unless a certain special situation occurs, 
and in \secref{T} we show that the homotopy can be constructed
even in that situation. By taking a product of all such cyclic groups
we obtain an embedding of $T(\ll)$ into a finite abelian group
$G(+)$. In \secref{4} this observation is put into the context
of a general theory describing group modifications of partial quasigroups
\cite{dk1}. We transfer the theory from noncommutative to abelian groups
and note that up to isomorphism there exists a unique finite abelian group 
$G(+)$ with a natural embedding $T(\ll) \to G(+)$ such that
any homotopy $T(\ll) \to H(+)$ can be factorized over this natural
embedding. Most of the needed facts are reproved, and that makes
the paper practically self-contained.

\section{Concepts and definitions}\label{2}

We have already mentioned that a dissection of an equilateral triangle 
is called \emph{separated} if
\begin{enumerate}
\item[(i)] no vertex of the dissection is of valence six, and
\item[(ii)] every line $p$ induced by a side of a dissecting triangle
has the property that the subset $U$ of $p$ forms a contiguous
segment if $U$ is defined as the union of all sides of dissecting
triangles that induce the line $p$.
\end{enumerate}

We shall sometimes consider the mates $T^\ll$ and $T^\rr$ of a 
latin bitrade $T$ as two partial quasigroups and write
$a_1\ll a_2 = a_3$ when $(a_1,a_2,a_3) \in T^\ll$ and 
$b_1 \rr\, b_2 = b_3$ when $(b_1,b_2,b_3) \in T^\rr$.

A triple of mappings $(\sigma_1,\sigma_2,\sigma_3)$ will be called a
\emph{homotopy} of partial quasigroups $T(\ll)$ and $S(\cdot)$ if
$\sigma_3(u_1\ll u_2) = \sigma_1(u_1)\cdot \sigma_2(u_2)$ whenever
$u_1\ll u_2$ is defined. Say that $T(\ll)$ can be \emph{embedded}
into $S(\cdot)$ if there exists a homotopy $(\sigma_1,\sigma_2,\sigma_3)$
such that all three mappings $\sigma_i$ are injective.

Call a latin bitrade $T=(T^\ll,T^\rr)$ \emph{indecomposable} if
it cannot be expressed as a disjoint union of two nonempty latin bitrades.
We shall investigate only indecomposable latin bitrades.

The number of triples in $T^\ll$ is called the \emph{size} of $T$.
Thus $s= |T^\ll| = |T^\rr|$. Latin bitrades are often represented
by tables as partial latin squares, and so $a_1$ can be regarded as
(a label of) a row, $a_2$ as a column and $a_3$ as a symbol, for 
every $(a_1,a_2,a_3) \in T^\ll$. Denote by $m$ the aggregate number
of rows, columns and symbols. Thus $m = o_1+o_2+o_3$, where 
$o_i = |\{\alpha;$ $\alpha = a_i$ for some $(a_1,a_2,a_3) \in T^\ll\}|$.

An (indecomposable) latin bitrade is said to be \emph{spherical}
if $m = s+2$. This equality can be derived from the Euler identity
and expresses the fact that the bitrade can be represented upon a sphere.
In this paper the topological aspects of latin bitrades will be 
needed only in \secref{T} where one can find more information.
Nevertheless,
we remark here that those indecomposable latin bitrades that yield an oriented 
surface in 
a direct way are called \emph{separated}. These are exactly those
bitrades in which there does not exist a row (a column, or a symbol)
such that one could construct a new bitrade of the same size by
dividing all cells of the row (the column, or the symbol) 
into two new rows (columns, symbols). Spherical bitrades are always
separated and the inequality $m\le s+2$ holds for every (indecomposable)
latin bitrade.

With each latin bitrade $T= (T^\ll, T^\rr)$ associate a set of 
equations $\eq(T)$ in such a way that every triple $(a_1,a_2,a_3) \in T^\ll$
yields the equation $a_1+a_2 = a_3$. The elements $a_i$ are thus
regarded as unknowns. It is natural to have different unknowns
for rows, columns and symbols, and so we assume that $a_i \ne b_j$
whenever $(a_1,a_2,a_3),(b_1,b_2,b_3) \in T^\ll$ and $1 \le i < j \le 3$.
(If the condition is violated,  then $T$  can be replaced by
an isotopic bitrade for which it is satisfied.)

To build a dissection one needs a latin
bitrade $T = (T^\ll, T^\rr)$ and a triple $a = (a_1,a_2,a_3) \in T^\ll$.
The pair $(T,a)$ will be called a \emph{pointed} bitrade.
In this section we shall always assume that $T$ is spherical.
We start from the set of equations $\eq(T,a)$ that is obtained
from $\eq(T)$ by removing the equation $a_1+a_2 = a_3$, and
by substituting $a_1=0$, $a_2=0$ and $a_3 =1$. We get in this way
a system with $m-3 = s-1$ variables and $s-1$ equations. 
In \secref{4} we shall explain why this system has always a
unique solution in rational numbers. If the solution to $b_i$ is $u_i$,
where $(b_1,b_2,b_3) \in T^\ll$, $b\ne a$,
write $\bar b_i=u_i$. In particular,
$\bar a_1 = \bar a_2 = 0$ and $\bar a_3 = 1$. Since $b_1$ is
usually associated with a row and $b_2$ with a column, we associate
in the euclidean plane the triple $b$ with the point 
$\bar b = (\bar b_2,\bar b_1)$.
We shall write $\tilde b_1$ to denote the line $y = \bar b_1$, and further
$\tilde b_2$ for the line $x = \bar b_2$ and $\tilde b_3$ for the
line $x+y = \bar b_3$. Note that $\tilde a_1$ is the axial line $y =0$,
$\tilde a_2$ is the axis $x=0$ and $\tilde a_3$ is the line 
$x+y = 1$. Denote by $\Sigma$ the triangle induced by these three
lines. Thus $\Sigma = \{(x,y);$ $x\ge 0$, $y\ge 0$ and $x+y \le 1
\}$. For each $c= (c_1,c_2,c_3)\in T^\rr$ denote by $\Delta = \Delta(c,a)$
the triangle induced by lines $\tilde c_1$, $\tilde c_2$ and $\tilde
c_3$. Of course, it is not clear that $\Delta$ is really
a triangle, i.~e.~that the lines $\tilde c_i$ do not meet in a single
point. If this happens, then we shall say that $\Delta$ \emph{degenerates}.
Our intention is to consider the question of when 
the triangles $\Delta(c,a)$, $c \in T^\rr$, dissect the triangle $\Sigma$
and none of them degenerates. 
In such a dissection each triangle has a side parallel to the axis
$x = 0$ and a side parallel to the axis $y = 0$. To obtain a dissection
of an equilateral triangle apply the affine transformation
$(x,y) \mapsto (y/2 + x, \sqrt 3 y/2)$.

Say that the solution to $\eq(T,a)$ is \emph{separated} if no two
unknowns representing two rows (or two columns, or two symbols) attain
the same value. In other words we require that $\bar b_i = \bar d_i$
if and only if $b_i = d_i$, for every $(b_1,b_2,b_3), (d_1,d_2,d_3) 
\in T^\ll$
and every $i \in \{1,2,3\}$. It is clear that if the solution to 
$\eq(T,a)$ is separated, then $\Delta(c,a)$ degenerates
for no $c\in T^\rr$.

Let us return to the procedure described in the beginning of \secref{1},
assuming that there is no dissecting vertex of valence six. The dissection
of $\Sigma$ can be then interpreted as a pointed bitrade $(T,a)$, 
where $a$ is the triple of lines that induce the sides of $\Sigma$.
Denote by $s$ the number of dissecting triangles and by $\ell$
the number of contiguous segments. We shall also use $m$ 
to denote  the aggregate number of rows, columns and symbols.
We see that $s$ is equal to the size of $T$, that $m \le \ell$, and 
that the
equality $m=\ell$ holds if and only if the dissection is separated.
Furthermore, there are $s+2$ dissecting vertices, and each of them
is an extreme point of exactly two edges. Thus $\ell = s+2$, and
$T$ is spherical if and only if the dissection is separated.

We can assume that $\Sigma =\{(x,y);$ $x\ge 0$, $y\ge 0$ and $x+y
\le 1\}$. The dissection yields a solution to $\eq(T,a)$. Since
there is only one solution, we see that the dissection can be
reconstructed as $\{\Delta(c,a);$ $c\in T^\rr\}$. A separated
dissection thus induces a separated solution to $\eq(T,a)$. 
\tref{21} claims that this observation can be reversed.
We shall prove it in \secref{3}. (Nearly all proofs in Sections 
\ref{2}--\ref{T} depend upon the fact that the linear system $\eq(T,a)$ 
always has a unique solution. That was proved already in \cite{dk3} and 
we prove it anew in \secref{4}, in \lref{43}.)

\begin{thm}\label{21}
Let $T=(T^\ll, T^\rr)$ be a spherical latin bitrade,
and suppose that $a=(a_1,a_2,a_3) \in T^\ll$ is such a triple
that the solution to $\eq(T,a)$ is separated. Then the set
of all triangles $\Delta(c,a)$, $c \in T^\rr$, dissects the triangle
$\Sigma = \{(x,y);$ $x\ge 0$, $y\ge 0$ and $x+y \le 1\}$. This
dissection is separated.
\end{thm}

The construction of latin bitrades from triangle
dissections was first described in \cite{d1}. It is also discussed
in \cite{kee}, and \cite{d2} gives a topological interpretation.
By \cite{d2} one can derive a spherical latin bitrade from dissections with
values of valence six as well. That opens a possibility to 
extend \tref{21} to all dissections.

For an example of a dissection, consider the following spherical bitrade 
$(T^{\ll},\, T^{\rr})$:  \\ 
\[
T^{\ll} = 
\begin{array}{|c||c|c|c|c|c|}
\hline \ll & c_0 & c_1 & c_2 & c_3 & c_4\\
\hline \hline r_0 & s_4 & ~ & s_0 & ~ & s_2\\
\hline r_1 & ~ & ~ & ~ & s_2 & s_4\\
\hline r_2 & s_0 & s_1 & s_2 & s_3 & ~\\
\hline r_3 & s_1 & s_3 & ~ & s_4 & ~\\
\hline \end{array} 
\quad 
T^{\rr} = 
\begin{array}{|c||c|c|c|c|c|}
\hline \rr & c_0 & c_1 & c_2 & c_3 & c_4\\
\hline \hline r_0 & s_0 & ~ & s_2 & ~ & s_4\\
\hline r_1 & ~ & ~ & ~ & s_4 & s_2\\
\hline r_2 & s_1 & s_3 & s_0 & s_2 & ~\\
\hline r_3 & s_4 & s_1 & ~ & s_3 & ~\\
\hline \end{array}
\label{eqnsphericalbitrade}
\]\\

\begin{figure}[htbp]
\begin{center}
 
\includegraphics[scale=0.75]{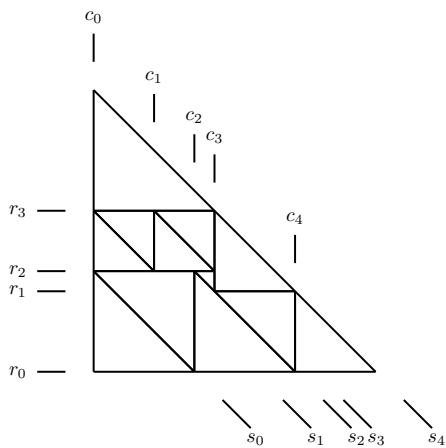}
  
\caption{Dissection for a spherical bitrade. The labels $r_i$, $c_j$,
$s_k$, refer to lines $y = \bar r_i$,
$x = \bar c_j$,
$x+y = \bar s_k$, respectively. }
\label{fg2}
\end{center}
\end{figure}

Let $a = (a_1,a_2,a_3) = (r_0,c_0,s_4)$. Then the system of equations
$\textnormal{Eq}(T, a)$ has the solution
\begin{align*}
\bar r_0 &= 0,\, \bar r_1 = 2/7,\, \bar r_2 = 5/14,\, \bar r_3 = 4/7 \\
\bar c_0 &= 0,\, \bar c_1 = 3/14,\, \bar c_2 = 5/14,\, \bar c_3 = 3/7,\, \bar c_4 = 5/7 \\
\bar s_0 &= 5/14,\, \bar s_1 = 4/7,\, \bar s_2 = 5/7,\, \bar s_3 = 11/14,\, \bar s_4 = 1.
\end{align*}
The dissection is shown in Figure~\ref{fg2}.
 Entries of
$T^{\rr}$ correspond to triangles in the dissection. For example,
$(r_0,c_0,s_0) \in T^{\rr}$ is the triangle bounded by the lines
$y = \bar r_0 = 0$, 
$x = \bar c_0 = 0$, 
$x+y = \bar s_0 = 5/14$ while 
$(r_1,c_3,s_2) \in T^{\ll}$ corresponds to the intersection of the
lines
$y = \bar r_1 = 2/7$, 
$x = \bar c_3 = 3/7$, 
$x+y = \bar s_2 = 5/7$.

We shall conclude this section by three preparatory lemmas. 
All of them assume that $T = (T^\ll, T^\rr)$ is a spherical
latin bitrade and that $a = (a_1,a_2,a_3)$ is a fixed element
of $T^\ll$.

By a \emph{lateral} point of a triangle we understand a point that
is upon a side of the triangle, but not equal to a vertex. Each 
triangle is thus a disjoint union of its vertices, lateral points and 
interior points.

In Lemmas~\ref{22} and~\ref{24} we shall use permutations 
$\mu_{r,s}$ of $T^\rr$, for every $r,s \in \{1,2,3\}$, $r \ne s$.
 Each of the mappings $\mu_{r,s}$
permutes the set $T^\rr$ in such a way that for every $c = (c_1,c_2,c_3)\in 
T^\rr$ we take first $d=(d_1,d_2,d_3) \in T^\ll$ such that $d_t = c_t$
for $t \ne s$, and then we choose $c'=\mu_{r,s}(c)\in T^\rr$ so that 
$c' = (c'_1,c'_2,c'_3)$ satisfies $c'_t = d_t$ for $t \ne r$. 
Note that by reversing the process we get $c = \mu_{s,r}(c')$. Hence
$\mu_{s,r} = \mu_{r,s}\m$.
Note that if  $t \in \{1,2,3\}$
is different from both $r$ and $s$, then $c'_t = c_t$
and $\mu_{s,r}(c) = \mu_{s,t}\mu_{t,r}(c)$. 

\begin{lem}\label{22} 
Assume that $b= (b_1,b_2,b_3) \in T^\ll$, and that
$1 \le i < j \le 3$. If $\bar b$ is not a vertex of $\Sigma$, 
then there exists $c = (c_1,c_2,c_3) \in T^\rr$
such that $\Delta(c,a)$ does not degenerate, $\tilde b_i = \tilde
c_i$ and $\tilde b_j = \tilde c_j$.
\end{lem}
\dem 
We shall say that $c \in T^\rr$ \emph{degenerates at} $b$ 
if $\Delta(c,a)$ degenerates
and the lines $\tilde c_1$, $\tilde c_2$ and $\tilde c_3$ meet in the point
$\bar b$. 
Denote by $C$ the set of all such $c \in T^\rr$. If 
$C = \emptyset$, then the lemma is obvious. Let us have $C$ nonempty,
and suppose first that there exist $r,s \in \{1,2,3\}$, $r\ne s$, such that
$\mu_{r,s}(c) \notin C$ for some $c \in C$.

While investigating this case we shall assume that the triples $d$,
$c$ and $c'$ have the same meaning as in the definition of $\mu_{r,s}$
and that $\{1,2,3\} = \{r,s,t\}$.
Since $c$ degenerates at $b$, there must be $\bar b = \bar d$. We
assume that $\bar b$ is not a vertex of $\Sigma$, and hence
$d \ne a$. Thus
$\tilde b_t = \tilde d_t = \tilde c'_t$ and $\tilde b_s = \tilde d_s
= \tilde c'_s$, by the construction of $c'$. Since $\Delta(c',a)$
does not degenerate, we are done when $\{i,j\}=\{s,t\}$. 
By switching $r$ and $s$ we also obtain immediately the case
$\{i,j\} = \{r,t\}$ since from $\mu_{s,r} =\mu_{r,s}\m$ we see
that there exists $e\in C$ such that $\mu_{s,r}(e) \notin C$.
This means that we also have a proof for the case $\{i,j\} = \{r,s\}$ 
when there exists $e\in C$ such that 
$\mu_{t,s}(e) \notin C$ or $\mu_{r,t}(e) \notin C$. However, we can always 
put $e = c$ or $e = \mu_{t,s}(c)$ since 
$\mu_{r,s}(c) = \mu_{r,t}\mu_{t,s}(c)\notin C$.
 
Suppose now that $C$ is nonempty and that $\mu_{r,s}(c) \in C$ for
every $c \in C$ and $r,s \in \{1,2,3\}$, $r \ne s$. Define $D$
as the set of all $d=(d_1,d_2,d_3) \in T^\ll$ such that $d$
agrees in  two coordinates with some $c \in C$. 
Choose $d \in D$ and $r \in \{1,2,3\}$. Let $c'\in T^\rr$
agree with $d$ in the two coordinates different from the $r$th
coordinate. We shall prove that $c'\in C$. There exists some
$c \in C$ that agrees with $d$ in two coordinates. If $c' = c$,
then there is nothing to prove. We can thus assume that 
$c_t=c'_t = d_t$, $c'_s = d_s \ne c_s$ and $c_r = d_r\ne c'_r$, 
where $\{r,s,t\} = \{1,2,3\}$.
But then $c' = \mu_{r,s}(c)$ and we have $c'\in C$ by our assumption.
We see that $(D,C)$ is a latin bitrade. Every $c \in C$ degenerates
at $b$. Since $\bar b$ is not a vertex of $\varSigma$, there cannot
be $a \in D$. Hence $(C,D)$  does not
coincide with $T=(T^\ll, T^\rr)$. This is a contradiction
because $T$ is assumed to be indecomposable. 
\edem

\begin{lem} \label{23} 
The triangle $\Delta(c,a)$ is contained in the triangle $\Sigma$, 
for every $c = (c_1,c_2,c_3) \in T^\rr$.
\end{lem}
\dem 
Put $H_1 = \{(x,y);$ $y<0\}$, $H_2 = \{(x,y);$ $x<0\}$ and
$H_3 = \{(x,y);$ $x+y > 1\}$. Our task is to obtain a contradiction
when $\bar b \in H_i$ for some $b =(b_1,b_2,b_3)\in T^\ll$, where
$b \ne a$ and $i \in \{1,2,3\}$. The three cases are similar and so
we can assume that there exists $b$ with $\bar b \in H_3$. 
Put $h = \max\{\bar b_3;$ 
$\bar b \in H_3\}$ and choose $b\in T^\ll$ such that
$\bar b_3 = h$ and $\bar b_1$ attains the maximum 
possible value. By \lref{22} there exists $c=(c_1,c_2,c_3)\in T^\rr$ such that 
$\Delta =\Delta(c,a)$ does
not degenerate, $\bar c_3 = h = \bar b_3$ and $\bar c_1 = \bar b_1$. 
Let $d,e\in T^\ll$ be such that $d = (d_1,c_2,c_3)$ and $e = (c_1,c_2,e_3)$,
for some $d_1$ and $e_3$. Then $\bar b$, $\bar d$ and $\bar e$
are the (pairwise distinct) vertices of $\Delta$. The vertex
$\bar d$ is upon the line $\tilde c_3 = \tilde b_3$, and hence
$h = \bar d_1 + \bar c_2 < \bar b_1 + \bar c_2$, by the maximality
of $\bar b_1$. The vertex $\bar e$ is upon the line $\tilde b_1 = 
\tilde c_1$, and hence $h > \bar c_1 + \bar c_2 = \bar b_1 + \bar c_2$.
We have obtained a contradiction.
\edem

\begin{lem} \label{24} 
Suppose that $b = (b_1,b_2,b_3) \in T^\ll$ and $i \in \{1,2,3\}$
are such that $b_i \ne a_i$. If $\bar b_i = \bar a_i$, then
$\Delta(c,a)$ degenerates whenever $c = (c_1,c_2,c_3) \in T^\rr$
and $c_i = b_i$. The point $\bar b$ is never a vertex of $\Sigma$.
\end{lem}
\dem 
Let us consider the case $i = 1$; the other cases are similar. Thus
$b_1 \ne a_1$ and $\bar b_1 = \bar a_1 = 0$. 
If $b' = (b_1,b_2',b_3') \in T^\ll$, then $\bar b_3' = \bar b_1
+ \bar b_2' = \bar b_2'$ since $b' \ne a$. Put $H = \{\bar b_2';
$ $(b_1,b_2',b_3') \in T^\ll\}$ and choose $h \in H$. 

Then there certainly exists at least one $c = (b_1,c_2,c_3) \in T^\rr$
such that $\bar c_2 = h$. Every  $c'= (b_1,c_2',c_3') \in T^\rr$
can be expressed as $\mu_{3,2}^i(c)$ for some $i \ge 0$ since 
the bitrade $T$ is separated. Can $c$ be chosen in such a way 
that $\Delta(c,a)$ does not degenerate?  If $\Delta(c,a)$ degenerates
and $c' = \mu_{3,2}(c)$, then $\bar c_2'= \bar c_3  = \bar c_2 = h$ as
well since $(b_1,c_2',c_3) \in T^\ll$. Therefore either $H = \{h\}$,
or the sought choice is possible. 

Let us thus suppose that $\Delta(c,a)$ does not degenerate.
Then $\bar c_2 \ne \bar c_3$ and $H$ contains at least two elements. 
Assume that $h = \max H$ and that $c'$ equals $\mu_{3,2}(c)$.
Both $(\bar c_2,0) = (h,0)$ and $(\bar c_2',0)$ are vertices 
of $\Delta(c,a)$. The third vertex is equal to $(\bar c_2,\bar d_1)$,
where $d=(d_1,c_2,c_3) \in T^\ll$. There cannot be $d =a$ since
$\bar c_2 = h > 0$.  From \lref{23} we get $\bar d_1 > 0$, and
hence $h=\bar c_2' = \bar c_3 = \bar d_1 + \bar c_2 > \bar c_2 = h$,
which contradicts the choice of $h$. We have proved that $H$
always contains only one element. That means that $\Delta(c,a)$
always degenerates when $c = (c_1,c_2,c_3) \in T^\rr$, $c_1 \ne a_1$
and $\bar c_1 = 0$ (more generally, when $c_i \ne a_i$ and $\bar c_i 
= \bar a_i$). 

What remains is to exclude cases $H =\{0\}$ and $H = \{1\}$.
Each of them
implies that the set $K_j = \{(u_1,u_2,u_3) \in T^\ll;$ $
|\bar u_j - \bar a_j| = 1\}$ is nonempty for some $j \in \{1,2,3\}$.
We shall show that this never occurs.

We shall
first do so under an additional assumption that the set
$C = \{c\in T^\rr;$ $\Delta(c,a)$ does not degenerate and
$\Delta(c,a) \ne \Sigma\}$ is nonempty. Put 
$\Gamma = \bigcup(\Delta(c,a);$ $c \in C)$ and set
$h_1 = \min \{y;$ $(x,y) \in \Gamma\}$ and $h_2 = \min\{x;$ 
$(x,h_1) \in \Gamma\}$. The point $(h_2,h_1)$ must be a vertex
of some $\Delta(e,a)$, $e \in C$. Hence there exists $b=(b_1,b_2,b_3)
\in T^\ll$ with $\bar b_1 = h_1$, $\bar b_2 = h_2$, and $\bar b_3 =
h_1+h_2$. Assume $(h_2,h_1) \ne (0,0)$. From \lref{23} we see that
$(h_2,h_1)$ is not a vertex of $\Sigma$, and from \lref{22} we
obtain the existence of $c = (c_1,c_2,c_3) \in C$ with 
$\bar c_1 = h_1$ and $\bar c_3 = h_1 + h_2$. Any triangle with
sides upon $\tilde c_1$ and $\tilde c_3$ contains an element $(x,y)$
such that either $y<h_1$, or $y =h_1$ and $x<h_2$. Hence
$(h_1,h_2) = (0,0)$ and there exists $e = (e_1,e_2,e_3) \in C$
with $\bar e_1 = \bar e_2 = 0$. By the first part of the proof,
$\Delta(e,a)$ degenerates if $e_1 \ne a_1$ or $e_2 \ne a_2$. However,
$\Delta(e,a)$ does not degenerate as $e \in C$. Hence
$e = (a_1,a_2,e_3)$. 

We shall now use the fact that $\Delta(e,a)$ does not degenerate
to prove that $K_3 = \emptyset$. Assume $K_3 \ne \emptyset$ and note
that $K_3 = \{(u_1,u_2,u_3)\in T^\ll;$ $\bar u_3 = 0\}$. To get a 
contradiction we shall use the indecomposability of $T$. To prove
that $K_3$ induces a decomposition of $T$ into two latin bitrades
it suffices to verify that $\bar u_3'=0$ whenever $u = (u_1,u_2,u_3)
\in K_3$, $v = (v_1,v_2,v_3) \in T^\rr$ and $u' = (u_1',u_2',u_3')
\in T^\ll$ are such that both pairs $\{u,v\}$ and $\{v,u'\}$ 
disagree in exactly one coordinate. There is nothing to prove if
$u_3'= u_3$. The first two coordinates can be treated similarly,
and so we can assume that $u_1 = u_1' = v_1$. 
Thus $u' = \mu_{3,2}(u)$ or $u' = \mu_{2,3}(u)$. It suffices to
consider only the former case since $\mu_{2,3}(u)$ can be obtained
by iterative applications of $\mu_{3,2}$ which start at $u$ as 
$\mu_{2,3} = \mu_{3,2}\m$. 

Let us have $u' = \mu_{3,2}(u)$. Then 
$u_3 = v_3$, $u_1 = v_1 = u_1'$ and $v_2 = u_2'$. Two vertices
of $\Delta (v,a)$ are upon the line $\tilde v_3 = \tilde u_3$. The
points of the line satisfy $x+y=0$. Both vertices are thus equal to $(0,0)$,
by \lref{23}, and hence $\bar u'$, the third vertex of $\Delta(v,a)$,
is equal to $(0,0)$ as well. Therefore either $\bar u'_3 = 0$, or $u' = a$. 
If $u' =a$, then $v_1 = a_1$ and $v_2 = a_2$. Thus $v$ coincides with 
the triple $e$ of the preceding paragraph. However, we cannot have
$e = v$ since $\Delta(v,a)$ degenerates while $\Delta(e,a)$ does not.
Hence $\bar u'_3 = 0$ in all cases.

Finally, let us suppose that the set $C$ is empty. 
If $\Delta(c,a) = \Sigma$ for some $c = (c_1,c_2,c_3) \in T^\rr$,
then $\bar c_i = \bar a_i$ for all $i \in \{1,2,3\}$. 
In at least one case there must be $c_i \ne a_i$ since $c \ne a$.
Hence $\Delta(c,a)$ has to degenerate, by the first part of the proof.
That is a contradiction, and therefore
what remains is to solve the case when $\Delta(c,a)$ degenerates
for every $c \in (c_1,c_2,c_3) \in T^\rr$. Every such $c$
thus satisfies $\bar c_1 + \bar c_2 = \bar c_3$. If $c = 
(a_1,c_2,c_3) \in T^\rr$, $d = (a_1,c_2,d_3) \in T^\ll$ and
$c' = (a_1,d_2,d_3) \in T^\rr$, then $c_2 \ne a_2$ implies
$\bar d_2=\bar d_3=\bar c_2=\bar c_3$ since $\bar a_1 = 0$. 
The arrows $c \to c'=\mu_{2,3}(c)$
thus retain the value $\bar c_2 = \bar c_3$ for all
$(a_1,c_2,c_3) \in T^\rr$, with a possible exception of one such $c$.
However, the arrows form a cycle and so there is no exception. 
By setting $d = a$ we see that there exists a cycle of $\mu_{2,3}$
in which both $c_2 = a_2$ and $c_3 = a_3$ occur. That supplies a contradiction
since $\bar a_2 \ne \bar a_3$.
Hence there always exists
$c \in T^\rr$ that does not degenerate.
\edem

\section{Separated solutions yield dissections}\label{3}

Suppose that the conditions of \tref{21} are satisfied.
The triple $a = (a_1,a_2,a_3)\in T^\ll$ is fixed throughout the
section and so we shall write $\Delta(c,a)$ simply as $\Delta(c)$, for
every $c \in T^\rr$.
Our goal is to show that 
these triangles dissect the triangle $\Sigma$. Denote the
vertices of $\Sigma$ by $A = (0,0)$, $B = (0,1)$ and $C= (1,0)$. 
By \lref{23},  $\Delta \subseteq \Sigma$ for all $\Delta(c)$,
$c \in T^\rr$.

In the rest of this section we shall call elements of $\{\Delta(c);$
$c \in T^\rr\}$
just \emph{triangles}. None of them degenerates because the solution
to $\eq(T,a)$ is assumed to be separated.
The union of all triangle sides will be denoted by $\Pi$,
and we shall be investigating the set $\Gamma = (\operatorname
{Int} \Sigma) - \Pi$.
For each $X \in \Gamma$ denote by
$\pi(X)$ the number of triangles in which $X$ is an interior point.
The set $\Gamma$ is open and we clearly have $\pi(X) = \pi(Y)$
when $X$ and $Y$ are in the same component of $\Gamma$.
We wish to show that the triangles dissect $\Sigma$, and that
is the same as proving $\pi(X) = 1$ for each $X \in \Gamma$.
We shall assume that 
$$\hat \Gamma = \{X \in \Gamma;\ \pi(X) \ne 1\} \ne \emptyset,$$
and argue by contradiction. The set $\hat \Gamma$ is open as well,
and its closure is compact. Therefore there exists a unique point
$P=(\alpha,\beta)$ in the closure of $\hat \Gamma$
such that if $P' = (\alpha',\beta')$ is another
point of the closure, then either (1) $\alpha+\beta < \alpha'+\beta'$,
or (2) $\alpha + \beta = \alpha'+\beta'$ and $\beta<\beta'$.

If $P$ is a lateral point of $\Sigma$, then $P$ is clearly incident
to one of the two axes. Furthermore, in such a case $P$ has to be a vertex
of a triangle since otherwise we might move $P$ towards $A$.

For every $X\in \Gamma$ there exists a neighbourhood $U \subseteq \Gamma$ 
in which $\pi$ is a constant function. Hence $P \notin \Gamma$.
It follows that $P\in \Pi$ since either $P \in (\operatorname{Int}
\Sigma) \setminus \Gamma$, or $P \in \Sigma \setminus
(\operatorname{Int}\Sigma)$. 

\begin{lem} \label{31}
The point $P$ does not equal $A$, $B$ or $C$.
\end{lem}
\dem Clearly $P\notin \{B,C\}$, by the definition of $P$. Assume
$P=A$ and consider the triple $a' = (a_1,a_2,a_3')\in T^\rr$
(recall that $\tilde a_1$ is the horizontal
axis and $\tilde a_2$ is the vertical axis). All points $X \in \Gamma$
that are close enough to $A$ belong to $\Delta(a')$ and that means 
that $\pi(X) \ge 1$. We assume $A = P$, and so there
must exist another triangle for which $A$ is a vertex. However,
that is impossible since such a triangle would have to be of the
form $\Delta(a'')$, where $a'' = (a_1,a_2,a_3'')\in T^\rr$ and $a_3' \ne
a_3''$.
\edem

Denote by $p_1$, $p_2$ and $p_3$ the lines that pass through $P$
and are parallel to $\tilde a_1$, $\tilde a_2$ and $\tilde a_3$, respectively.
By \lref{31} the existence of a triple $b = (b_1,b_2,b_3)\in T^\ll$
such that $\tilde b_i = p_i$ for all $i \in \{1,2,3\}$ is equivalent 
to the fact that $P$ is a triangle vertex.

The next lemma is an immediate consequence of the assumption that 
the solution to $\eq(T,a)$ is separated. 

\begin{lem} \label{32}
Assume that $P$ is a triangle vertex. Then for each $i$ and $j$,
where $1\le i < j \le 3$, there exists a unique $c = (c_1,c_2,c_3)
\in T^\rr$ such that $\tilde c_i = p_i$ and $\tilde c_j = p_j$.
We shall denote $\Delta(c)$ by $\Delta(i,j)$.
\end{lem}

\begin{lem} \label{33}
The point $P$ is not a lateral point of $\Sigma$.
\end{lem}
\dem Assume the contrary. Since $P$ cannot be upon $\tilde a_3$, it has to
be an axial point, as we have already remarked in the text above.
Denote by $p$ the axis that passes through $P$ and by $p^+$ the half-plane
induced by $p$ that contains $\Sigma$. If $P$ were not a triangle vertex,
we could move $P$ along $p$ towards $A$. Hence $P$ is a triangle vertex.
By \lref{32} for all $i$ and $j$ with $1 \le i < j \le 3$ there 
exists a triangle $\Delta(i,j)$. Clearly $\Delta(i,j) \subset p^+$.
We can thus find a neighbourhood $U$ of $P$ such that 
$U\cap \Sigma \subset \Delta(1,2) \cup \Delta(1,3) \cup \Delta(2,3)$.
No two of these three triangles share an interior point and no
other triangle has $P$ as a vertex. Since $P$ is in the closure of 
$\hat \Gamma$, there must exist a triangle $\Delta$ that intersects
every neighbourhood of $P$. We see that $P$ is a lateral point of $\Delta$,
and therefore there exists a neighbourhood $U$ of $P$ such that
$U\cap \Delta \cap \Delta(i,j) \ne \emptyset$ for all choices of $i$ and $j$.
That makes possible a move of $P$ within $U$ along $p$ towards $A$,
and that is our concluding contradiction.
\edem

\begin{figure}[htbp]
\begin{center}
 
\input{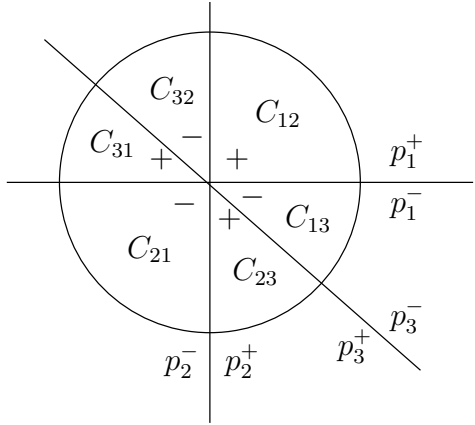}
  
\caption{The sectors and half-planes induced by the point $P$.}
\label{fg1}
\end{center}
\end{figure}

We can thus assume that $P$ is an interior point of $\Sigma$. Choose a circular
neighbourhood of $P$ and consider the six sectors of the circle
that are induced by lines $p_i$, $1 \le i \le 3$. 
Call the sectors alternately \emph{positive} and \emph{negative}
so that we regard as positive the sector formed by all points
$(\alpha',\beta')$ with $\alpha' \ge \alpha$ and $\beta' \ge \beta$.
Of the two sectors bordered by
$p_i$ and $p_{i+1}$ denote the positive by $C_{i,i+1}$  
and the negative by $C_{i+1,i}$. The indices are computed modulo 3
and the neighbourhood
is chosen small enough to have $\operatorname{Int}C_{ij}\subseteq \Gamma$. 
Thus $\pi(X) = \pi(X')$ when there exist $i,j \in \{1,2,3\}$
such that both $X$ and $X'$ belong to $\operatorname{Int} C_{ij}$, 
and we shall denote this integer by $\pi_{ij}$. The arrangement
of sectors is illustrated by Figure~\ref{fg1}.

For $i \in \{1,2,3\}$ denote by $p_i^+$ the half-plane determined by $p_i$ 
that contains $C_{i,i+1} \cup C_{i-1,i+1} \cup C_{i-1,i}$, and by
$p_i^{-}$ the opposite half-plane.

\begin{lem}\label{34}
If $(i,j) \in \{(3,1), (2,1),
(2,3), (1,3)\}$, then $\pi_{ij} = 1$.
\end{lem}
\dem Suppose first that $\pi_{ij} \ge 2$ or $\pi_{ij} = 0$, where
$(i,j) \in \{(3,1), (2,1), (2,3)\}$. Then 
$\operatorname{Int}C_{ij}\subseteq \hat \Gamma$ and
$\alpha' + \beta' < \alpha + \beta$ for every 
$(\alpha',\beta') \in \operatorname{Int}C_{ij}$.
If $\pi_{13} \ge 2$ or $\pi_{13} = 0$, then $\alpha' + \beta' = \alpha + \beta$
and $\beta' < \beta$ for every point $(\alpha',\beta')\in C_{13} \cap p_3$
that differs from $P$. In every
case we thus obtain a contradiction to the choice of $P$.
\edem

There may exist a triangle $\Delta$ such that $C_{ij} \subset \Delta$
and $P$ is a vertex of $\Delta$. If such a triangle exists, then it is
determined uniquely, by \lref{32}, and we shall denote it by $\Delta_{ij}$. 
If the triangle exists, then it is equal to $\Delta(i,j)$ when $i < j$,
and to $\Delta(j,i)$ when $j< i$. Note that
if $\Delta_{ij}$ exists, then $\Delta_{ji}$ does not exist.

By writing $\Delta \subset p_i^{\pm}$ we mean that either
$\Delta \subset p_i^{+}$, or $\Delta \subset p_i^{-}$.
If the point $P$ is a lateral point of a triangle $\Delta \subset
p_i^{\pm}$, then $C_{rs} \subset \Delta$ for every $C_{rs} \subset p_i^\pm$.
Furthermore, $\pi_{rs} = 1$ for at least one such $(r,s)$, by \lref{34}.
This means that $\Delta$, if it exists, is uniquely determined,
and we shall denote it by $\Delta_i^{\pm}$ (the existence of $\Delta_i^+$ need
not exclude the existence of $\Delta_i^-$).

For formal reasons it is sometimes useful to write $p_i^{\veps}$ in
place of $p_i^+$ or $p_i^-$, where $\veps \in \{-1,1\}$. Under this notation
$$p_i^\veps \supset C_{i,i+\veps} \cup C_{i-\veps,i+\veps} \cup 
C_{i-\veps,i} \text{\, for every \,} \veps \in \{-1,1\} \text{\, and \,}
i \in \{1,2,3\}.$$ 

\begin{lem}\label{35}
Let $i \in \{1,2,3\}$ and $\veps \in \{-1,1\}$ be such that 
both $\Delta_{i,i+\veps}$ and $\Delta_{i-\veps,i}$ exist and
neither of them is equal to $\Delta_{32}$ or $\Delta_{12}$. 
If $\pi_{i-\veps,i+\veps}\ge 1$, then the triangle $\Delta_{i-\veps,
i+\veps}$ exists as well.
\end{lem}
\dem 
Let $\Delta$ be a triangle that contains $C_{i-\veps,i+\veps}$, and
assume that  $\Delta_{i-\veps,i+\veps}$ does not exist. Then $P$
has to be a lateral or interior point of $\Delta$, and so 
$\Delta$ contains $C_{i,i+\veps}$ or $C_{i-\veps,i}$. Hence
$\pi_{i,i+\veps} \ge 2$ or $\pi_{i-\veps,i} \ge 2$, and that
contradicts our assumptions, by \lref{34}.
\edem

\begin{lem}\label{36}
Suppose that $P$ is a vertex of a triangle. Then for some
$i \in \{1,2,3\}$ and $\veps \in \{-1,1\}$ there exist triangles
$\Delta_{i,i+\veps}$, $\Delta_{i-\veps,i+\veps}$, $\Delta_{i-\veps,i}$
and $\Delta_i^{-\veps}$.
\end{lem}
\dem By \lref{32} there exist exactly three triangles of the 
form $\Delta_{jk}$. Suppose that they are not contained in any
half-plane $p_i^{\veps}$. Then they have to correspond either
to all sectors $C_{j,j+1}$ (the positive sectors), or to all sectors
$C_{j,j-1}$ (the negative sectors). However, in both these cases
we easily obtain a contradiction by means of Lemmas~\ref{35} and~\ref{34}.
Therefore there exists a unique pair $(i,\veps)$ such that all
three triangles with vertex $P$ are contained in $p_i^\veps$.
Note that the sectors $C_{i,i+\veps}\cup C_{i,i-\veps}$ and
$C_{i+\veps,i} \cup C_{i-\veps,i}$ are opposite, and 
that they are separated by interspersed sectors $C_{i-\veps,i+\veps}$
and $C_{i+\veps,i-\veps}$. Hence no union of two adjacent sectors $C_{rs}$
intersects both $C_{i,i+\veps}\cup C_{i,i-\veps}$ and
$C_{i+\veps,i} \cup C_{i-\veps,i}$. It follows that 
one of them has no interior point common with the sector $C_{12}\cup
C_{32}$. Thus  $\pi_{i,i+\veps} = \pi_{i,i-\veps} = 1$ or
$\pi_{i+\veps,i} = \pi_{i-\veps,i} = 1$, by \lref{34}. 
Let $\Delta$ be the unique triangle
that contains $C_{i,i-\veps}$ in the former case, and $C_{i+\veps,i}$
in the latter case. The triangle $\Delta$ has been chosen in such a way
that it is not contained in $p_i^\veps$. All three triangles for
which $P$ is a vertex are contained in $p_i^\veps$. Hence $P$ is
not a vertex of $\Delta$ and has to be a lateral or interior point
of $\Delta$. We shall choose a triangle $\Delta' \subset p_i^\veps$
for which $P$ is a vertex and which overlaps with no other triangle
upon a neighbourhood of $P$. In the former case put 
$\Delta' = \Delta_{i,i+\veps}$ (and use $\pi_{i,i+\veps}=1$), while
in the latter case set $\Delta' = \Delta_{i-\veps,i}$ (and use
$\pi_{i-\veps,i} = 1$). Triangles $\Delta$ and $\Delta'$ have been
defined in such a way that $\Delta \cap \Delta' \cap p_i$ is a nontrivial
segment. The choice of $\Delta'$ guarantees that they do not overlap,
and therefore $p_i$ has to induce a side of $\Delta$. This means
that $P$ is a lateral point of $\Delta$, and so $\Delta = \Delta_i^{-\veps}$.
\edem

\begin{lem}\label{37}
The point $P$ is a vertex of no triangle. 
\end{lem}
\dem
Let $i$ and $\veps$ be as in \lref{36}. The point $P$ clearly cannot be
a vertex or an interior point of any triangle not listed in that lemma.
If there would exist $\Delta_h^\eta \ne \Delta_i^{-\veps}$, then
necessarily $\pi_{rs} \ge 2$ for some $(r,s) \in  \{(3,1), (2,1),
(2,3), (1,3)\}$, which is a contradiction to \lref{34}. 
Therefore $P$ cannot be
a lateral point either.
\edem

To finish the proof of \tref{21} it therefore suffices to
consider the case when $P$ is a lateral point of one or more
triangles, but a vertex of no triangle. All triangles in which $P$
is a lateral point are of the form $\Delta_h^{\pm}$,
$1 \le h \le 3$. From \lref{34} we see that there exist at 
least two such triangles and that $P$ is not an interior point
of any triangle. Suppose that among the triangles there exists
a pair $(\Delta_h^+, \Delta_h^-)$. Note that such a pair can be found always
when there are at least four triangles. If there were no other triangles
beyond the pair,
then every sector $C_{ij}$ would be covered by exactly one triangle.
If there were further triangles, then at least one of the sectors 
listed in \lref{34} would be covered by at least two triangles.
Hence no such pair exists.

Suppose now that $P$ is contained in exactly two triangles. By
\lref{34} they have to intersect in one of the sectors $C_{12}$ and $C_{32}$.
But then the opposite sector (i.~e.~$C_{21}$ or $C_{23}$) is covered by 
no triangle, and we obtain again a contradiction to \lref{34}.


We see that $P$ must be included in
exactly three triangles. We have nine pairs $(C_{ij},
\Delta_h^\veps)$ with $C_{ij}\subset \Delta_h^\veps$. By \lref{34}, only 
$C_{32}$ and $C_{12}$ can appear in more than one pair. 
Since $9 > 6 + 1 +1$, at least one of them
has to appear exactly three times.  But then the opposite
sector lacks covering, and we obtain another contradiction
with \lref{34}. We have proved that the point $P$ does not exist,
and therefore $\hat \Gamma=\emptyset$.

\section{Trigons} \label{T}

In \lref{T2} we shall observe that every solution to $\eq(T,a)$ can
be interpreted, for some $n \ge 2$,  as a homotopy of $T(\ll)$ 
to $\mathbb Z_n(+)$. Our goal is to show that this set of homotopies
is rich enough to separate any two elements within the same group
(the three groups are the rows, the columns and the symbols). This 
is nearly achieved by \lref{24}. To solve the remaining cases
we have to consider certain configurations within
latin bitrades which we shall call
\emph{trigons}. We shall first describe them informally.

Suppose that a bitrade $S_0$ is derived from a dissection of $\Sigma$ and that
$\Delta$ is a dissecting triangle. Let us consider another dissection,
say of $\Sigma'$, that yields a bitrade $S_1$. Use the latter
dissection as a pattern for how to dissect $\Delta$. By identifying
vertices of $\Sigma'$ with vertices of $\Delta$ we obtain in this way
a new dissection of $\Sigma$.
The dissection determines a new bitrade, say $T$, that can
be regarded as a superimposition of $S_0$ and $S_1$. 
The ensuing definition stipulates that the triple representing
$\Delta$ in $S_0^\rr$ becomes a trigon in $T$ (while in $S_0$ it is not
regarded as a trigon).

A \emph{trigon} in a latin bitrade $T=(T^\ll,T^\rr)$ is 
a triple $c = (c_1,c_2,c_3)$, $c \notin T^\ll$, for which there
exist elements $c_i'\ne c_i$ such that $(c_1,c_2,c_3')$, $(c_1,c_2',c_3)$
and $(c_1',c_2,c_3)$ belong to $T^\rr$. These three triples will
be called the \emph{corner triples} of $c$.

In \lref{T3} we shall show how to recombine a homotopy 
$\vhi: S_0 \to \mathbb Z_n$ and a homotopy $\psi: S_1 \to \mathbb
Z_m$ into a homotopy $T \to \mathbb Z_{nm}$. If $S_0$
and $S_1$ are the bitrades based on dissections that have been described above,
then one can say that the new homotopy is obtained from
$\vhi$ by embedding $\mathbb Z_n$ into $\mathbb Z_{nm}$ and refining
the images of points within $\Delta$ by means of $\psi$. This
construction will provide an inductive tool for how to separate
an element within a trigon from an element outside. In 
\lref{T4} we shall show that trigons always arise in those
cases when \lref{24} does not suffice to separate $\bar b_i$ and $\bar a_i$. 
This makes the inductive argument possible, but first we
have to show that any trigon $c$ induces trades
$S_0$ and $S_1$ that are smaller than $T$. In other words,
we have to explain how the above construction of $T$ from $S_0$
and $S_1$ can be repeated when $T$ is not defined by means of dissections.

That is not difficult but it requires a digression that explains how
to treat spherical bitrades by means of the combinatorial topology
\cite{d2,kee,ldcd}.
The standard way \cite{cav,cl}
is to associate with a separated latin bitrade $T=(T^\ll,T^\rr)$
a black and white triangulation in which every $(a_1,a_2,a_3) \in
T^\ll$ is turned into a white triangle $\{a_1,a_2,a_3\}$ and
every $(b_1,b_2,b_3)\in T^\rr$ yields a black triangle $\{b_1,b_2,b_3\}$.
Here we shall use an alternative way how to associate with $T$ a
combinatorial surface. The alternative might be called the
\emph{semidual} of the black and white triangulation. It
can be described either directly \cite{d2,d3}, or by the ensuing 
modification of the standard definition.  (Note that if the latin bitrade
$T$ is not separated, then the above definition of the black and
white triangles yields only a pseudosurface. The kissing points
of such a pseudosurface correspond to elements $a_t$, $t \in
\{1,2,3\},$ that can be divided into two or more rows (columns, symbols)
in such a way that the new structure is also a latin bitrade.)

The points of the semidual are the elements of $T^\ll$. They represent
the centers of the white triangles. There are two kinds of 
faces---the \emph{cyclic} faces and the \emph{triangular} faces.
Each cyclic face is induced by a (unique) element $a_t$ that
represents a row, a column or a symbol, and is formed by the
cyclic sequence of white triangle centers (which we identify
with the elements of $T^\ll$) of those triangles that have $a_t$
as a vertex. The number  
of cyclic faces thus equals the aggregate number of rows, columns
and symbols. The triangular faces correspond in a one-to-one manner
to black faces (and thus also to elements of $T^\rr$). They are 
formed by unordered triples that consist of the centers
of the three white triangles that are adjacent to the given
black triangle. In other words every $(c_1,c_2,c_3) \in T^\rr$
induces a triangular face formed by the uniquely determined
points $(c_1',c_2,c_3)$, $(c_1,c_2',c_3)$, $(c_1,c_2,c'_3) \in T^\ll$.

We shall now define permutations $\nu_{r,s}$ of $T^\ll$ that are dual
to the already defined permutations $\mu_{r,s}$ of $T^\rr$. The cyclic
faces are exactly the cycles (in the cyclic decomposition) of $\nu_{r,s}$.
Assume that $r,s \in \{1,2,3\}$ and that $r\ne s$. For $a =(a_1,a_2,a_3)
\in T^\ll$ consider $(b_1,b_2,b_3)\in T^\rr$ and $a' = (a_1',a_2',a_3')
\in T^\ll$ such that $a_i \ne b_i$ exactly when $i = s$ and $a_i' \ne 
b_i$ exactly when $i = r$. Then $a' = \nu_{r,s}(a)$. Choose $t$
so that $\{r,s,t\} = \{1,2,3\}$. It is clear that the cycle
of $\nu_{r,s}$ that passes through $a$ does not change  $a_t$
and coincides with the cyclic face induced by $a_t$. 

We have $\nu_{r,s}\nu_{s,t} = \nu_{r,t}$ and $\nu_{r,s}\m = \nu_{s,r}$.
For every $j \in \{1,2,3\}$ denote $\nu_{j+1,j-1}$ by $\tau_j$ (the
indices are computed modulo 3). Note that $\tau_1\tau_2\tau_3 =
\tau_2\tau_3\tau_1 = \tau_3\tau_1\tau_2$ is the identity mapping
and that the permutations $\tau_j$ induce an orientation of cyclic
faces that can be used to orient the combinatorial surface.

Let now $T = (T^\ll,T^\rr)$ be a separated latin bitrade and let
$c = (c_1,c_2,c_3)$ be a trigon in $T$. Let $\alpha_j \in T^\ll$
be the triple that agrees with $c$ in coordinates $j\pm 1$ and
let $\gamma_j \in T^\rr$ be the corner triples of $c$ that
also agree with $c$ in these coordinates. 

When we follow the definition of $\nu_{j\pm 1,j}(\alpha_j)$,
we see that $\alpha_j$ is first changed to $\gamma_j$, and then
a change is made in the $(j\pm 1)$th coordinate. Hence $\nu_{j\pm 1,j}
(\alpha_j)$ are the vertices, together with $\alpha_j$, of
the triangular face that is associated with $\gamma_j$. 
They can be expressed as $\nu_{j-1,j}(\alpha_j) = \tau_{j+1}(\alpha_j)$
and $\nu_{j+1,j}(\alpha_j) = \tau\m_{j-1}(\alpha_j)$.
In particular, $\tau_j(\alpha_{j-1})$ and
$\tau_j\m(\alpha_{j+1})$ are vertices of the triangular faces
that are induced by $\gamma_{j-1}$ and $\gamma_{j+1}$, respectively.

Let $\ell_j$ be the length of the cyclic face that is induced by $c_j$.
Both $\alpha_{j-1}$ and $\alpha_{j+1}$ are incident to the cycle,
and hence $\tau_j^{k_j}(\alpha_{j-1}) = \alpha_{j+1}$ for some 
positive $k_j < \ell_j$. We cannot have $\tau_j(\alpha_{j-1}) = \alpha_{j+1}$
as $c \ne \gamma_{j-1}$, and so $k_j \ge 2$. This gives us the
following characterization of trigons (the converse implication
is clear):

\begin{lem}\label{T1}
Let $T = (T^\ll,T^\rr)$ be a separated latin bitrade. A triple
$c = (c_1,c_2,c_3)$ is a trigon in $T$ if and only if for
all $j \in \{1,2,3\}$ there exist $\alpha_j \in T^\ll$ that
differ from $c$ exactly in the $j$th coordinate, and there
exist integers $k_j$ such that
$$\tau_j^{k_j}(\alpha_{j-1}) = \alpha_{j+1},\ \ 2 \le k_j < \ell_j,$$
where $\ell_j$ is the length of the cycle of $\tau_j$ that moves
$\alpha_{j\pm 1}$.
\end{lem}

\lref{T1} does not require that $T$ were spherical. However, if it is not,
then the oriented closed path (it will be denoted by $P$)
$$ \alpha_3,\tau_1(\alpha_3)\dots,\tau_1^{k_1}(\alpha_3) = \alpha_2,
\tau_3(\alpha_2),\dots,
\tau_3^{k_3}(\alpha_2) = \alpha_1,\tau_2(\alpha_1),
\dots,\tau_2^{k_2}(\alpha_1) =\alpha_3$$
need not separate the surface into two disjoint parts.
Since we are interested in spherical latin bitrades we can 
assume that the separation takes place.

Put $\beta_j=\tau_{j-1}^{k_{j-1}-1}(\alpha_{j+1})$ and consider
the subpath $\beta_j,\alpha_j,\tau_{j+1}(\alpha_j)$ of $P$. We shall
denote it by $P_j$,
$j \in \{1,2,3\}$. We have $\tau\m_j(\beta_j) = \tau_{j+1}(\alpha_j)$  
since
$$\tau_j(\tau_{j+1}(\alpha_j)) = \tau_j\tau_{j+1}\tau_{j-1}(\beta_j) =
\beta_j.$$
Let $h_j$ be the length of the cycle of $\tau_j$ that moves $\beta_j$. 
Then $$ \beta_j,\tau_j(\beta_j),\dots,\tau_j^{h_j-1}(\beta_j) =
\tau_{j+1}(\alpha_j)$$ is an oriented path from $\beta_j$ to
$\tau_{j+1}(\alpha_j)$ that will be denoted by $Q_j$. 
By substituting $Q_j$ for $P_j$ in $P$, 
$j \in \{1,2,3\}$, we obtain an oriented closed path that is called
the \emph{inner circumference} of the trigon $c$.

The points upon the circumference and inside are the \emph{inner points}
of $c$. The other elements of $T^\ll$ are the \emph{outer points} of $c$.
The outer points include the (vertex) points $\alpha_1$, $\alpha_2$ and 
$\alpha_3$. 

The triangular faces inside the inner circumference determine a subset
of $T^\rr$ and we define $S_1^\rr$ as the union of this subset with
$\{\gamma_1,\gamma_2,\gamma_3\}$. As $S_1^\ll$ take the union of
$\{c\}$ with the set of all inner points. Then $S_1 = (S_1^\ll,S_1^\rr)$
forms a latin bitrade, and we shall call it the \emph{inner bitrade} of the
trigon $c$.

Similarly, define 
$S_0^\ll$ as the set of all outer points and $S_0^\rr$
as the set obtained by unifying $\{c\}$ with the set of all elements in
$T^\rr\setminus\{\gamma_1,\gamma_2,\gamma_3\}$ that determine a
triangular face outside the inner circumference. Note that all
vertices of every such face are outer points. Hence 
$S_0 = (S_0^\ll,S_0^\rr)$ forms a bitrade as well, and we shall call
it the \emph{outer bitrade} of the trigon $c$.

To understand the meaning of $S_0$ and $S_1$ topologically is easy. 
Let $T$ be a spherical bitrade. The trigon $c$ determines upon the
combinatorial sphere of the semidual a triangular area that is 
described by the
path $P$. The sphere of $S_0$ is obtained by deleting the inner
structure of this area which is now considered as a new triangular
face. The sphere of $S_1$ is obtained by deleting everything outside the area
and merging $\alpha_1$, $\alpha_2$ and $\alpha_3$ into a single
new point. This merge converts the three sides of the triangular
area into cyclic faces. 

We have seen that if $T$ is spherical, then both $S_0$ and $S_1$
are spherical as well. Of course, this can be also proved
formally \cite{d4} by counting the faces and using the eulerian 
characteristic. 

This finishes our digression and we return to the program outlined
in the beginning of this section. We shall use the conventions
of \secref{2} established for a situation in which there are fixed
an indecomposable spherical
latin bitrade $T = (T^\ll,T^\rr)$ and a triple $a=(a_1,a_2,a_3)
\in T^\ll$. 
(For example, if $b = (b_1,b_2,b_3) \in T^\ll$,
then $\bar b_j$ denotes the value assigned to $b_j$ by
$\eq(T,a)$ for every $j \in \{1,2,3\}$, and $\bar b =
(\bar b_2,\bar b_1)$.)

\begin{lem}\label{T2}
Let $T= (T^\ll,T^\rr)$ be a spherical latin bitrade and let
$a=(a_1,a_2,a_3)$ be an element of $T^\ll$. Then there exists
an integer $n\ge 2$ and a homotopy $\psi = (\psi_1,\psi_2,\psi_3)$
of $T(\ll)$ into $\mathbb Z_n$ such that for any
$(b_1,b_2,b_3),(d_1,d_2,d_3)\in T^\ll$ and $j \in \{1,2,3\}$
the equality $\psi_j(b_j) = \psi_j(d_j)$ holds if and only if 
$\bar b_j = \bar d_j$ (i.e.~if there coincide the values
assigned to $b_j$ and $d_j$ by the linear system $\eq(T,a)$).
\end{lem}
\dem Let $n$ be the least positive integer such that $n\bar b_j
\in \mathbb Z$ for every $b = (b_1,b_2,b_3) \in T^\ll \setminus \{a\}$
and every $j \in \{1,2,3\}$. Then $n\bar b_1 + n \bar b_2 =
n \bar b_3$, $n\bar a_1 + n \bar a_2 = 0$ and $n\bar a_3 = n$. 
Hence by defining $\psi = (\psi_1,\psi_2,\psi_3): T^\ll \to 
\mathbb Z_n$ so that $\psi_j(d_j) \equiv n\bar d_j \bmod n$
for every $(d_1,d_2,d_3) \in T^\ll$ and every $j \in \{1,2,3\}$
we indeed obtain a homotopy. Note that $0 \le n\bar d_j < n$
if $j \in \{1,2\}$ and $0 < n\bar d_3 \le n$, by Lemmas~\ref{23}
and~\ref{24}. Hence $n\bar b_j \equiv n \bar d_j \bmod n$ implies
$\bar b_j = \bar d_j$. By \lref{24} none of $\Delta(c,a)$ equals 
$\Sigma$, and so $n\ge 2$.
\edem

The homotopy $\psi$ described in the proof of \lref{T2} will be
called the homotopy \emph{induced} by $a$. Denote the integer
$n\bar d_j$ by $\bar \psi_j(d_j)$, for every $d=(d_1,d_2,d_3) \in 
T^\ll$ and every $j \in \{1,2,3\}$. Then
$$ \bar \psi_1(d_1) + \bar \psi_2(d_2) = \bar \psi_3(d_3)
\text{\, whenever \,} d \ne a.$$
The triple $(\bar \psi_1,\bar \psi_2,\bar \psi_3)$ will be called
the \emph{near-homotopy} induced by $a$. We shall also say that
$n$ is the \emph{width} of $a$ in $T$.

\begin{lem}\label{T3}
Let $T=(T^\ll,T^\rr)$ be a spherical latin bitrade with a trigon
$c = (c_1,c_2,c_3)$. Let $S_0$ and $S_1$ be the outer and inner
bitrades of $c$, let $\vhi = (\vhi_1,\vhi_2,\vhi_3)$ be a homotopy
$S_0(\ll) \to \mathbb Z_m$, let $n$ be the width of $c$ in $S_1$
and let $(\bar\psi_1,\bar\psi_2,\bar\psi_3)$ be the near-homotopy
induced by $c \in S_1^\ll$. Denote by $\rho$ the embedding
$\mathbb Z_m \to \mathbb Z_{mn}$, $i \mapsto ni$, put
$h_1 = \rho\vhi_1(c_1)$, $h_2 = \rho \vhi_2(c_2)$ and $h_3 =
h_1+h_2$. Finally, choose $k \in \mathbb Z_{mn}$ so that 
$$ k \equiv \vhi_3(c_3) - \vhi_1(c_1) - \vhi_2(c_2) \bmod m.$$
For $j \in \{1,2,3\}$ and $b = (b_1,b_2,b_3) \in T^\ll$ put
$\vhi'_j(b_j) = \rho\vhi_j(b_j)$ if $b$ is an outer point
of $c$ and $\vhi'_j(b_j) =h_j + \bar \psi_j(b_j)k$ if $b$
is an inner point of $c$. Then $\vhi' = (\vhi'_1,\vhi'_2,\vhi'_3)$
is a homotopy $T(\ll) \to \mathbb Z_{mn}$. 
\end{lem}
\dem The two formulas that define $\vhi'_j$ overlap in $c_j$.
The case $j \in \{1,2\}$ is disambiguous since then $\bar \psi_j(c_j)
= 0$. From $h_3=h_1+h_2$ and $\bar \psi_3(c_3)=n$ we obtain that 
$$h_3 + \bar \psi_3(c_3) k \equiv n(\vhi_1(c_1) + \vhi_2(c_2) + k)
\bmod mn.$$
Since $k \equiv \vhi_3(c_3) - \vhi_1(c_1) - \vhi_2(c_2) \bmod m$
we see that $h_3 + \bar \psi_3(c_3)k \equiv n\vhi_3(c_3)\bmod mn$,
and so $h_3 + \bar \psi_3(c_3)k = \rho \vhi_3(c_3)$.

The definition of $\vhi'$ is hence correct and $\vhi'_1(b_1) + \vhi'_2(b_2)$
clearly equals $\vhi'_3(b_3)$ if $b = (b_1,b_2,b_3) \in T^\ll$ is an outer
point. If $b$ is an inner point, then $\bar \psi_1(b_1) + \bar \psi_2(b_2)
= \bar \psi_3(b_3)$ and so the equality holds as well.
\edem

Fix an (indecomposable) spherical latin bitrade $T=(T^\ll,T^\rr)$ and
a triple $a = (a_1,a_2,a_3) \in T^\ll$. 
Denote by $a'_j$, $j \in \{1,2,3\}$, the elements such that
the triples $\eta_1 = (a_1',a_2,a_3)$, $\eta_2 = (a_1,a'_2,a_3)$
and $\eta_3 = (a_1,a_2,a_3')$ belong to $T^\rr$. Note that
$\mu_{j-1,j+1}(\eta_{j+1})=\eta_{j-1}$ and denote by $k_j$
the length of the respective cycle. Then $\mu^{k_j - 1}_{j-1,j+1}
(\eta_{j-1}) = \eta_{j+1}$ and $k_j \ge 2$. 

In the proof of \lref{T4} some phrases will be expressed in 
a shortened way. By saying that $c = (c_1,c_2,c_3)
\in T^\rr$ \emph{degenerates at} $(s,t)\in \Sigma$ we shall mean
that $\Delta(c,a)$ degenerates and $(s,t) = (\bar c_2,\bar c_1)$. 
We shall be also saying that an element $b_j$ \emph{shrinks
at} $(s,t)$ if at $(s,t)$ there degenerates every 
$c = (c_1,c_2,c_3) \in T^\rr$ with $c_j = b_j$. By \lref{24},
if $b_j \ne a_j$ and $\bar b_j = \bar a_j$, then $b_j$ shrinks
at some lateral point of $\Sigma$.

\begin{lem}\label{T4}
The triangle $\Delta(\eta_j,a)$ degenerates for no $j \in \{1,2,3\}$.
Fix $j \in \{1,2,3\}$ and define integers $0 \le i_0 < \dots 
< i_\ell = k_j - 1$ as those $i$, $0 \le i < k_j$, for which
$\Delta(\mu_{j-1,j+1}^i(\eta_{j-1}),a)$ does not degenerate.
Put $\gamma_r = \mu_{j-1,j+1}^{i_r}(\eta_{j-1})$, $0 \le r \le \ell$,
and denote by $b_r^\pm$ the $(j\pm 1)$th coordinate of $\gamma_r$.
For $0 \le r < \ell$ define $\beta_r = (f_1,f_2,f_3)$ so that $f_j =
a_j$, $f_{j-1} = b_r^-$ and $f_{j+1} = b_{r+1}^+$. Then either
$\beta_r \in T^\ll$ and $i_{r+1} = i_r+1$, or $\beta_r$ is a trigon
and $i_{r+1} > i_r + 1$. If $f' = (f_1',f_2',f_3') \in T^\rr$ 
agrees with $\beta_r$ in two coordinates, then $\Delta(f',a)$ 
never degenerates.

If $b=(b_1,b_2,b_3) \in T^\ll$ and $j \in \{1,2,3\}$ are such that
$\bar b_j = \bar a_j$ and $b_j \ne a_j$, then there exists unique
$r$, $0\le r < \ell$, such that  $\beta_r$ is a trigon and $b$
is an outer point of $\beta_r$. 
\end{lem}
\dem If $b \in T^\ll$ and $b \ne a$, then $\bar b$ is never equal
to a vertex of $\Sigma$, by \lref{24}. Therefore no $c \in T^\rr$
degenerates at a vertex of $\Sigma$, and hence $\eta_j$ cannot
degenerate for any $j \in \{1,2,3\}$. 

Assume now that $j =1$; the other cases are similar. 
By the definition, $\gamma_r = (a_1, b_r^+,b_r^-)$ whenever $0 \le r
\le \ell$. Suppose that $r < \ell$ and define $c_2$ and $c_3$
by $(a_1,c_2,b_r^-) \in T^\ll$ and $(a_1,c_2,c_3)=\mu_{3,2}(\gamma_r)
\in T^\rr$. Then $(a_1,b^+_{r+1},b^-_r) \in T^\ll$ $\Leftrightarrow$
$c_2 = b_{r+1}^+$ $\Leftrightarrow$ $(a_1,b_{r+1}^+,c_3) \in T^\rr$
$\Leftrightarrow$ $\mu_{3,2}(\gamma_r) = \gamma_{r+1}$. Thus
$\beta_r \in T^\ll$ if and only if $i_{r+1} = i_r + 1$.

All triangles $\Delta(\gamma_r,a)$ are inside $\Sigma$, by \lref{23}. 
One side of such a triangle is upon the horizontal axis. Another side,
which is upon the left, is parallel to the vertical axis and the third
side is parallel to the line $x+y = 1$. The triangles $\Delta(\gamma_r,a)$
and $\Delta(\gamma_{r+1},a)$ have a common vertex, but that obviously
is the only nonempty intersection of $\Delta(\gamma_r,a)$ and 
$\Delta(\gamma_s,a)$ whenever $0 \le r < s \le \ell$. Hence
if $b = (b_1,b_2,b_3) \in T^\rr \setminus \{a\}$ is such that
$\bar b = (h,0)$ is upon the horizontal axis, then $\bar b$ 
has to coincide, by Lemmas~\ref{22} and~\ref{24}, with the 
common vertex of $\Delta(\gamma_r,a)$ and $\Delta(\gamma_{r+1},a)$,
for some $r$, $0 \le r < \ell$. The value of $r$ will be now regarded
as fixed. Thus $h = \bar b_{r}^- = \bar b_{r+1}^+ = \bar c_2$.

We shall need two auxiliary claims:

\vspace{-0.7em}
\subsection*{Claim A} Suppose that $u = (u_1,u_2,b_r^-) \in T^\rr$
is such that $u_1$ shrinks at $(h,0)$. Let $(u_1',u_2',b_r^-)$
be equal to $\mu_{1,2}(u)$. Then either $u_2' = b_{r+1}^+$, or both
$u_1'$ and $u'_2$ shrink at $(h,0)$. 

\smallskip
The proof of Claim A depends upon 

\vspace{-0.7em}
\subsection*{Claim B} Suppose that $v = (v_1,u_2',v_3) \in T^\rr$
is such that $v$ degenerates at $(h,0)$. Let $v' = (v_1',u_2',v_3')$
be equal to $\mu_{3,1}(v)$. Then either $v'$ degenerates at $(h,0)$,
or $v' = \gamma_{r+1}$.

\smallskip
The assumptions of Claim B imply $\bar v_1 = 0$, $\bar u_2' = \bar v_3
= h$ and $(v_1',u_2',v_3) \in T^\ll$. Thus $\bar v_1' = 0$,
and $v_1'$ shrinks at $(h,0)$ if $v_1' \ne a_1$, by \lref{24}. 
Assume $v_1' = a_1$ and suppose that $v'$ does not degenerate.
Then $v'$ has to be equal to one of $\gamma_s$, $0 \le s \le \ell$,
and from $\bar u_2' = h$ we see that there must be $s = r+1$. 

With Claim B proved we can turn to Claim A. Assume $u_2' \ne b^+_{r+1}$.
We have $(u_1,u_2',b_r^-)
\in T^\ll$, and thus $u' = (u_1,u'_2,u_3) \in T^\rr$ for some $u_3$.
Now, $(u'_1,u'_2,b_r^-) = \mu_{1,3}(u')$, and therefore 
$\mu_{1,2}(u) = \mu_{3,1}^{k'}(u')$, where $k'+1$ is the length
of the cycle of $\mu_{3,1}$ that passes through $u'$. The triple
$u'$ degenerates at $(h,0)$ since $u_1$ shrinks at $(h,0)$. 
No $\mu_{3,1}^i(u')$ equals $\gamma_{r+1}$ since $u_2' \ne b^+_{r+1}$. 
Hence repeated applications of Claim B imply
that $\mu_{3,1}^{k'}(u') = \mu_{1,2}(u)$ degenerates at $(h,0)$.
Thus $\bar u_1' = 0$, and $u'_1$ shrinks at $(h,0)$ if $u'_1
\ne a_1$, by \lref{24}. To finish the proof of Claim A it thus
suffices to show how to deduce a contradiction from the assumption
that $u_1' = a_1$. However, in such a case $\mu_{1,2} (u) = \gamma_r$. 
That cannot be since $\mu_{1,2}(u)$ degenerates while 
$\gamma_r$ does not.

\smallskip
Suppose now that $\beta_r \in T^\ll$. Then there exists a unique
$c_1 \ne a_1$ such that $c = (c_1,b_{r+1}^+,b_r^-) \in T^\rr$.
Note that $\gamma_r = \mu_{2,1}(c)$ and $\gamma_r = \mu_{1,2}^k(c)$,
where $k+1$ is the length of cycle of $\mu_{1,2}$ that passes through
$c$. If $1 \le i \le k$ and $(u'_1,u'_2,b_r^-) = \mu_{1,2}^i(c)$,
then $u_2'  \ne b_{r+1}^+$. If $c$ degenerates, then $c_1$ shrinks,
by \lref{24}, and repeated
applications of Claim A imply that $\gamma_r = \mu_{1,2}^k(c)$
degenerates as well. This is a contradiction, and hence $c$ 
does not degenerate.

Let us turn to the case $\beta_r \notin T^\ll$. Then $(a_1,c_2,c_3)
= \mu_{3,2}(\gamma_r)$ degenerates at $(h,0)$ and $c_2 \ne b_{r+1}^+$.
Repeated applications of Claim B imply that $c_2$ shrinks at $(h,0)$.
We have $(a_1,c_2,b_r^-) \in T^\ll$, and so there exists $c_1'$ such 
that $(c_1',c_2,b_r^-) \in T^\rr$ and $c_1' \ne a_1$. This triple
degenerates at $(h,0)$ since $c_2$ shrinks. \lref{24} implies
that $c_1'$ shrinks at $(h,0)$ as well. Apply now Claim A
repeatedly, starting from $(c_1',c_2,b_r^-)$. Since $b_r^-$ does not
shrink, we have to hit $b_{r+1}^+$ at some stage. Therefore there
exist $c_1''$ and $c_1$ such that $c_1''$ shrinks at $(h,0)$,
$(c_1'',b_{r+1}^+,b_r^-) \in T^\ll$, and 
$c = (c_1,b_{r+1}^+,b_r^-) \in T^\rr$. We see that $\beta_r$ is
a trigon with corner triples equal to $c$, $\gamma_r$ and $\gamma_{r+1}$.

Our next step is to show that $c$ does not degenerate. The proof is
practically the same as in the case $\beta_r\in T^\ll$. Assume that
$c$ degenerates and apply Claim A repeatedly. The permutation 
$\mu_{1,2}$ brings $c$ to $\gamma_r$ before one hits $b_{r+1}^+$,
and so $\gamma_r$ has to degenerate. This is a contradiction and hence
$c$ cannot degenerate.

To finish the proof of the lemma it thus remains to show that $b$
is an outside point of the trigon $\beta_r$ whenever $b = (b_1,b_2,b_3)
\in T^\ll$ is such that $\bar b = (h,0)$ and $b_1 \ne a_1$. We shall be thus
also proving that no such $b$ exists when $\beta_r \in T^\ll$.

Let $B^\ll \subseteq T^\ll$ be the set of all counterexamples (for a fixed
$h$) and
assume that $B^\ll \ne \emptyset$. Set $B = \{b_1;$ $(b_1,b_2,b_3)
\in B^\ll$ for some $b_2$ and $b_3 \}$ and note that if 
$b' = (b_1,b_2',b_3') \in T^\ll$ and $b_1 \in B$, then $\bar b'$ 
equals $(h,0)$
since $b_1$ shrinks at $(h,0)$ by \lref{24}. If $\beta_r$ is a trigon,
then for a given $b_1\ne a_1$ either all elements $b'$ are among the inside
points of the trigon or all such elements are outside points. 
Since we assume that $(b_1,b_2,b_3)$ is an inside point for some $b_2$
and $b_3$, the former alternative has to take place. It follows that 
$B^\ll = \{(b_1,b_2,b_3) \in T^\ll;$ $b_1 \in B\}$.
Put $B^\rr = \{(b_1,b_2,b_3)
\in T^\rr; $ $b_1 \in B\}$. We shall show that $(B^\ll,B^\rr)$
is a latin bitrade. This will yield the sought contradiction since
$T$ is assumed to be indecomposable.

Obviously is suffices to prove

\vspace{-0.5em}
\subsection*{Claim C} If $(e_1,e_2,e_3) \in B^\rr$ and 
$(d_1,e_2,e_3) \in T^\ll$, then $d_1 \in B$; and

\vspace{-0.5em}
\subsection*{Claim D} If $(d_1,d_2,d_3) \in B^\ll$ and
$(e_1,d_2,d_3) \in T^\rr$, then $e_1 \in B$.

\smallskip 
We shall first explain how Claim C implies Claim D. Suppose
that $(d_1,d_2,d_3) \in B^\ll$ and that $e = (e_1,d_2,d_3) \in 
T^\rr$. Let $d_3'$ be such that $d = (d_1,d_2,d_3') \in T^\rr$.
Then $e = \mu_{1,3}(d)$ and $d \in B^\rr$. By working inductively along
$\mu_{3,1} = \mu_{1,3}\m$ it is therefore enough to prove
that if $d' = (d_1',d_2,d_3'') \in B^\rr$, then $\mu_{3,1}(d')
= (d_1'',d_2,d_3''')$ belongs to $B^\rr$ as well. However,
this is clear from Claim C since $(d_1'',d_2,d_3'') \in T^\ll$. 

To prove Claim C put $e = (e_1,e_2,e_3)$ and $d= (d_1,e_2,e_3)$.
We have $e_1 \in B$, and hence $e_1$ shrinks at $(h,0)$, by 
\lref{24}. Therefore $e$ degenerates at $(h,0)$ and $\bar d$
equals $(h,0)$. 

If $\beta_r$ is a trigon, consider $d_2$ and $d_3$ such that
$(e_1,d_2,e_3) \in T^\ll$ and $(e_1,e_2,d_3) \in T^\ll$. 
Both these triples belong to $B^\ll$ and are inner points of $\beta_r$
since $e_1 \ne a_1$. The third triple from $T^\ll$ that is induced 
by $e = (e_1,e_2,e_3)$ is also an inner point of $\beta_r$, unless 
$e$ is the corner  triple. The third triple is equal to $d$,
of course. 
We have proved above that no corner triple
of $\beta_r$ degenerates, and therefore $d$ has to be an inner
point of $\beta_r$. The trigon $\beta_r$ has been constructed in such a way
that $\gamma_r$ and $\gamma_{r+1}$ are its corner triples and that
any $u=(a_1,u_2,u_3)\in T^\ll$ with $\bar u=(h,0)$ is an outer point
of $\beta_r$. Since $d$ is an inner point with $\bar d = (h,0)$,
there cannot be $d_1 = a_1$. Hence $d_1 \in B$.

If $\beta_r \in T^\ll$, then the argument is similar. From \lref{24} we
deduce that $d_1 \in B$ or that $d_1 = a_1$. If $d_1 = a_1$, then
$d=\beta_r$ since here we assume the existence of only one $u=(a_1,u_2,u_3)$ 
such that $\bar u = (h,0)$. If $d=\beta_r$, then $e$ is equal to the triple
$c = (c_1,b^+_{r+1},b^-_r)\in T^\rr$. We have already proved above 
that this triple does not degenerate. 

The proof is finished.
\edem

\begin{thm}\label{T5}
Let $T = (T^\ll,T^\rr)$ be a spherical latin bitrade.
If $a = (a_1,a_2,a_3)\in T^\ll$, $b = (b_1,b_2,b_3) \in 
T^\ll$ and $i \in \{1,2,3\}$ are such that $a_i \ne b_i$,
then there exist $n \ge 2$ and a homotopy $\vhi = (\vhi_1,
\vhi_2,\vhi_3): T(\ll) \to \mathbb Z_n(+)$ such that
$\vhi_i(a_i) \ne \vhi_i(b_i)$.
\end{thm}
\dem Consider first the homotopy $\psi = (\psi_1,\psi_2,\psi_3)$
induced by $a$. If $\psi_i(a_i) = \psi_i(b_i)$, then $\bar a_i
= \bar b_i$ by \lref{24}, and $b$ is an outer point of a trigon
$c = (c_1,c_2,c_3)$ in which $c_i = a_i$, by \lref{T4}.

In the rest we shall proceed by induction on the size of $T$.
Let $S_0$ be the outer trade of $c$, and let $S_1$ be the inner
trade. Recall that $c \in S_0^\rr$. By induction assumption
there exist $m \ge 2$ and a homotopy $\vhi = (\vhi_1,\vhi_2,\vhi_3):$
$S_0(\ll) \to \mathbb Z_m$ such that $\vhi_i(a_i) \ne \vhi_i(b_i)$.
\lref{T3} describes a procedure how to find $n'> n$ and a homotopy
$\vhi'=(\vhi'_1,\vhi'_2,\vhi'_3):$ $T(\ll) \to \mathbb Z_{n'}$
such that $\vhi'_i(a_i) \ne \vhi'_i(b_i)$.
\edem

\begin{cor}\label{newc}
Let $T = (T^\ll,T^\rr)$ be a spherical latin bitrade. Then both
$T(\ll)$ and $T(\rr)$ can be embedded into a finite abelian group.
\end{cor}
\dem By \tref{T5} for every triple $(a,b,j) \in T^\ll \times T^\ll \times
\{1,2,3\}$ such that $a_j \ne b_j$ there exists 
an integer $n = n[a,b,j]$ and a homotopy $\sigma=\sigma[a,b,j]: 
T(\ll) \to \mathbb Z_n$ in which $\sigma_j(a_j) \ne \sigma_j(b_j)$.
Put $G = \prod \mathbb Z_{n[a,b,j]}$ and define a homotopy $\tau:T(*) \to G$
so that the projection of $\tau$ to $\mathbb Z_{n[a,b,j]}$ coincides 
with $\sigma[a,b,j]$.
There are only finitely many triples $(a,b,j)$, and that makes $G$ finite.
The homotopy $\tau$ differentiates between any two different triples
of $T^\ll$ and hence it really embeds $T(\ll)$ into $G$.
\edem

%

\section{Modifications to abelian groups}\label{4}

Suppose that $T = (T^\ll,T^\rr)$ is a spherical
latin bitrade. $T$ is assumed to be indecomposable and 
we define $m = o_1 + o_2 + o_3$ in the same way as in 
\secref{2}. The linear system $\eq(T)$ has $m$ variables
and $s = m-2$ equalities. Relabel the variables by $x_1,\dots,x_m$
in such a way that the unknowns $x_j$ correspond to
rows, columns and symbols if $1 \le j \le o_1$, $o_1 < j \le o_1+o_2$
and $o_1 + o_2 < j \le m$, respectively. Order the equalities
of $\eq(T)$ in an arbitrary way and define a matrix
$B$ so that the $i$th row expresses the $i$th equation. (If
$x_r + x_s - x_t =0$ is the equation, then
$b_{ir}=b_{is}=1$, $b_{it} = -1$, and $b_{ij} = 0$ in other 
cases.)

Denote by $B_{ij}$ the matrix that is
derived from $B$ by omitting the $i$th and $j$th column.
By \cite[Lemma~3.3]{dk3} the following statement holds.
We include the proof since \cite{dk3} uses an approach that
is unnecessarily general for the needs of this paper.

\begin{lem}\label{41}
Suppose that either $1\le i \le o_1 < j \le m$ or $o_1 < i 
\le o_1+o_2 < j \le m$. Then $|\det B_{ij}| = |\det B_{1m}|$.
\end{lem}
\dem Let $C$ be a $(m+1)\times m$ matrix with rows $\mathbf c_1,
\dots, \mathbf c_{m+1}$ and let $\lambda_1,\dots,\lambda_{m+1}$
be coefficients such that $\sum \lambda_h\mathbf c_h = 0$. Denote
by $C_h$ the square matrix that is obtained from $C$ by deleting
the row $\mathbf c_h$. Consider $C_u$ and $C_v$ where $1\le u < v \le
m+1$, and denote by $D$ the matrix obtained from $C_v$ when the
$u$th row $\mathbf c_u$ is replaced by $-\lambda_v\mathbf c_v =
\sum_{h\ne v}\lambda_h\mathbf c_h$. Then $\det D = \lambda_u \det C_u
= -\lambda_v (-1)^{v-1-u}\det C_v$, and so $|\lambda_u \det C_u| =
|\lambda_v \det C_v|$.

Suppose now that $1\le r \le o_1 < s \le o_1+o_2 < t \le m$,
and form a $(m+1)\times m$ matrix $C$ from $B$ by adding $(m-1)$th, $m$th
and $(m+1)$th row such that each of them contains $m-1$ zeros, and the
cell in the column $r$ (or $s$, or $t$) contains the value $1$,
respectively. Finally, multiply the last row by $-1$.
Our intention is to show that $|\det C_u| = |\det C_v|$
when $m-1 \le u < v \le m+1$.

The determinants vanish when $B$ is not of rank $m-2$, and so we can
assume that it is of the full rank. Denote by $W$ the space of vectors
$(w_1,\dots, w_m)$ in which $w_1=\dots = w_{o_1}$, $w_{o_1+1} = \dots
= w_{o_1+o_2}$, $w_{o_1+o_2+1} = \dots = w_{m}$ and
$w_1 + w_{o_1+1} - w_{o_1+o_2+1} = 0$. Let $\mathbf w_1$
(or $\mathbf w_2$, or $\mathbf w_3$) be the element of $W$
such that $w_1 = 1$, $w_{o_1+1}=-1$ and $w_{o_1+o_2+1} = 0$
(or $w_1 = 0$ and $w_{o_1+1}=1 =w_{o_1+o_2+1}$, or
$w_{o_1+1} = 0$ and $w_{1}=1 =w_{o_1+o_2+1}$, respectively).
The rows $\mathbf c_1,\dots,
\mathbf c_{m-2}$ are orthogonal to the elements of $W$. Choose 
$\lambda_1,\dots,\lambda_{m+1}$ such that not all of them are zero
and $\sum \lambda_h \mathbf c_h = 0$. 
The scalar product of $\sum \lambda_h \mathbf c_h$ with any element
of $W$ thus vanishes, and hence $\mathbf u =\lambda_{m-1}
\mathbf c_{m-1} + \lambda_m \mathbf c_m + \lambda_{m+1}\mathbf c_{m+1}$
is also orthogonal to $W$. The scalar products of $\mathbf u$
with $\mathbf w_1$, $\mathbf w_2$ and $\mathbf w_3$ yield
$\lambda_{m-1} - \lambda_m
= \lambda_m - \lambda_{m+1} = \lambda_{m-1} - \lambda_{m+1} = 0$.
We cannot have $\lambda = \lambda_{m+1} = \lambda_m = \lambda_{m-1} = 0$
since $B$ is of rank $m-2$. Therefore $\lambda \ne 0$ and 
thus $|\det C_{m-1}| = |\det C_m| = |\det C_{m+1}|$, by the first
part of the proof.

Choose now $(r,s,t)$ in such a way
that $(i,j)$ is one of $(r,s)$, $(r,t)$ and $(s,t)$, and use the
obvious equalities
$|\det C_{m-1}| = |\det B_{st}|$, $|\det C_{m}| = |\det B_{rt}|$ and
$|\det C_{m+1}| = |\det B_{rs}|$. The rest is clear since we can move
in at most three steps from $(i,j)$ to any other $(i',j')$, including
the pair $(1,m)$.
\edem

We shall also need a result that follows from \cite[Lemma~3.1]{dk2}.
We give a full proof since \cite{dk2} seems to be difficult to 
read.  The proof is similar to the proofs of Sections~\ref{2}.
The homotopies described in \lref{42} will be called \emph{trivial}.

\begin{lem}\label{42}
Let $(\sigma_1,\sigma_2,\sigma_3)$ be a homotopy of $T(\ll)$
into the additive group of integers $\mathbb Z(+)$. Then
$\sigma_i(a_i) = \sigma_i(b_i)$ for
all $(a_1,a_2,a_3), (b_1,b_2,b_3) \in T^\ll$ and every 
$i \in \{1,2,3\}$.
\end{lem}
\dem Suppose that $(\sigma_1,\sigma_2,\sigma_3)$ is a nontrivial
homotopy. Say that $c = (c_1,c_2,c_3)\in T^\rr$ \emph{degenerates}
if $\sigma_1(c_1) + \sigma_2(c_2) = \sigma_3(c_3)$. Put
$M = \{(\sigma_1(a_1),\sigma_2(a_2),\sigma_3(a_3));$ $(a_1,a_2,a_3)
\in T^\ll\}$, $h_3 = \max \{r_3;$ $(r_1,r_2,r_3) \in M\}$,
$h_2 = \min \{r_2;$ $(r_1,r_2,h_3) \in M\}$, $h_1 = h_3 - h_2$,
$X = \{(a_1,a_2,a_3) \in T^\ll;$ $\sigma_i(a_i) = h_i$, 
$1 \le i \le 3\}$, and for each $i \in \{1,2,3\}$ define $Y_i$
as $\{(c_1,c_2,c_3)\in T^\rr;$ $\sigma_j(a_j) = h_j$ if
$1 \le j \le 3$ and $i\ne j\}$. 
For each element of $Y_i$ there exists a unique element of $X$
that agrees in two coordinates and disagrees in the $i$th
coordinate. An element $(c_1,c_2,c_3) \in Y_i$ degenerates
if and only if $\sigma_i(c_i) = h_i$, and therefore
$D = Y_1\cap Y_2 \cap
Y_3=Y_1\cap Y_2 = Y_1\cap Y_3 = Y_2 \cap Y_3$ consists of exactly 
those $c \in Y_1\cup Y_2 \cup Y_3$ that
degenerate. We have $|X| = |Y_i|$ for all $i$, and hence either 
$D = Y_1 = Y_2 = Y_3$, or $D$ is a proper subset
of each $Y_i$, $1 \le i \le 3$. In the former case all of the mappings 
$\mu_{r,s}$ permute $D$, and that makes $(X,D)$ a subtrade. An
indecomposable trade contains no proper subtrade, and hence
$(X,D) = (T^\ll, T^\rr)$. However, this is not possible since 
$(\sigma_1,\sigma_2,\sigma_3)$ is assumed to be nontrivial.

If $\mu_{r,s}(D) = D$, then there exists $c \in Y_i$ with $\mu_{r,s}(c)
\notin D$ since $D \subsetneq Y_i$. This inclusion guarantees the
existence of $c\in Y_i$ with $\mu_{r,s}(c) \notin D$ also in the
case when $\mu_{r,s} (D) \ne D$. 

Fix $c = (c_1,c_2,c_3) \in Y_1$ such that $\mu_{2,1}(c) =c'= (a_1,a_2,c_3)
\notin D$. Note that $(a_1,c_2,c_3) \in X$, $\sigma_1(a_1)
= h_1$, $\sigma_2(c_2) = h_2$ and $\sigma_3(c_3) = h_3$. Since
$c'$ does not degenerate, there must be $\sigma_2(a_2) \ne h_2$.

There exist elements $a_1'$ and $c'_3$ such that both
$(a_1,a_2,c_3')$ and $(a_1',a_2,c_3)$ belong to $T^\ll$. We have
$\sigma_3(c_3') \ne h_3$ since $\sigma_1(a_1) + \sigma_2(a_2) 
\ne h_1+h_2 = h_3$.
Therefore $\sigma_3(c_3')=h_1 + \sigma_2(a_2) < h_3$, 
and thus $\sigma_2(a_2)< h_2$.
On the other hand $\sigma_1(a'_1) + \sigma_2(a_2) = h_3$ implies
$\sigma_2(a_2) \ge h_2$, by the definition of $h_2$. We have obtained
a contradiction.
\edem


Let $G= G(+)$ be an abelian group. Consider, for a while, $\eq(T)$ as a
set of equations in $G$. We can regard each solution
as a triple $(\sigma_1,\sigma_2,\sigma_3)$ of mappings
into $G$ such that $\sigma_1(a_1) + \sigma_2(a_2) = \sigma_3(a_3)$
for all $(a_1,a_2,a_3) \in T^\ll$.
The solutions thus correspond to homotopies $T(\ll) \to G$.
Trivial homotopies can be obtained easily by choosing
elements $g$ and $h$ of $G$ and by setting $\sigma_1(a_1) = g$, 
$\sigma_2(a_2) = h$ and $\sigma_3(a_3) = g+h$, for all
$(a_1,a_2,a_3) \in T^\ll$. 

Consider $a = (a_1,a_2,a_3) \in T^\ll$ and assume 
that $a_1$, $a_2$ and $a_3$ have been relabelled as $x_r$,
$x_s$ and $x_t$, respectively. Thus $1 \le r \le o_1 <s
\le o_1+o_2 < t \le m$. If $\det B_{rs} = 0$, then there
exists a nonzero integer vector $v$ such that $Bv^\top = 0$ and
$v_r = v_s = v_t = 0$. Such a vector supplies a solution to $\eq(T)$
in $\mathbb Z$, and thus it yields a nontrivial
homotopy $T(\ll) \to \mathbb Z$. However, no such homotopy 
exists, by \lref{42}, and hence $\det B_{rs} \ne 0$. 
Remove now from $B$ the columns $r$, $s$, $t$ and the row
that corresponds to the equation $x_r + x_s = x_t$. Let it
be the $i$th row. The new matrix, say $C$, can be thus obtained 
from $B_{rs}$ by deleting a column and the $i$th row. Since
this row contains a single nonzero value, and since this value is 
equal to $\pm 1$ and is in the column that is being deleted we see that $|\det C| 
= |\det B_{rs}| \ne 0$.
The matrix $C$ is the matrix
of the linear system $\eq(T,a)$. 
We have hence proved the following statement:

\begin{lem}\label{43}
The system of linear equations $\eq(T,a)$ has a unique solution
in rational numbers. Furthermore, $\det B_{1m} \ne 0$.
\end{lem}


Note that results of Sections~\ref{2} and~\ref{T} assume
the validity of \lref{43}. Thus only at this point we can
regard as proved the fact that each spherical latin bitrade
can be embedded into a finite abelian group. We shall now
take a more systematic look upon such embeddings
and connect \cref{newc} to earlier results
of \cite{dk1} in which there was developed a machinery of group 
modifications. We shall limit our discussion only to aspects 
relevant to abelian groups, and refer to \cite{dk1} for 
further ramifications. The terminology of \cite{dk1} is somewhat 
different, but that should not cause difficulties.

We will view homotopies as morphisms in the category of partial
quasigroups. Both groups and latin bitrades can be regarded
as subcategories of this category: A group homomorphism
$f:G \to H$ is identified with a homotopy $(f,f,f)$, and a homotopy
$(\sigma_1,\sigma_2,\sigma_3)$ of $T(\ll) \to S(\ll)$ is regarded
(for our purposes here) as a morphism $(T^\ll,T^\rr) \to (S^\ll,S^\rr)$.



The following statement is clear and does not require a proof.
The notation $\sigma[a]$ that is introduced in \lref{45} will
be used further on (in particular, in and before \lref{47}).

\begin{lem}\label{45}
Let $T = (T^\ll,T^\rr)$ be a latin bitrade and let $\sigma = (\sigma_1,
\sigma_2, \sigma_3)$ be a homotopy of $T(\ll)$ into an abelian
group $K = K(+)$. For $a= (a_1,a_2,a_3)\in T^\ll$ define
a triple of mappings $\sigma[a] = (\tau_1,\tau_2,\tau_3)$ in such a way that
$\tau_i(b_i) = \sigma_i(b_i)-\sigma_i(a_i)$ for all $i \in \{1,2,3\}$
and all $(b_1,b_2,b_3) \in T^\ll$. Then $\sigma[a]$ is also
a homotopy of $T(\ll)$ into $K$ and $\tau_i(a_i) = 0$ for each
$i \in \{1,2,3\}$. The homotopy $\sigma$ is an embedding if and only
if $\sigma[a]$ is an embedding, and $\sigma$ is trivial if and only
if $\sigma[a]$ is trivial.
\end{lem}

By a \emph{modification} (or \emph{reflexion}) of a category into
a subcategory one understands morphisms $g_K: K \to G(K)$ that are
defined for each object $K$. The object $G(K)$ is in the subcategory
and for all morphisms $h: K \to H$ where the target $H$ is in the subcategory
there exists (in the subcategory) a unique morphism $k: G(K) \to H$
such that $h = kg_K$. As a typical examples one can take the
natural projections $G \to G/G'$ which yield a modification from
the category of groups into the subcategory of abelian subgroups.

Let $T=(T^\ll, T^\rr)$ be a latin bitrade. For every 
$(a_1,a_2,a_3)\in T^\ll$ regard
again elements $a_i$ as variables and denote
by $F = F(T)$ the free abelian group generated by them. The group
$F$ is thus of rank $m = o_1+o_2+o_3$, but we do not
require that the size of $T$ is necessarily equal to $m-2$. 
Regard now $\eq(T)$
as a set of elements $a_1+a_2-a_3 \in F$, denote by $N(T)$ the
subgroup generated by these elements, and denote by 
$\mathbf G(T)$ the factor-group $F(T)/N(T)$. Define $\mathbf g_T$
as a triple of mappings $(g_1,g_2,g_3)$ such that $g_i(a_i) = a_i + N(T)$.
Then $\mathbf g_T$ is clearly a homotopy of $T(\ll)$ into $\mathbf G(T)$,
and it can be verified easily that $\mathbf g_T: T \to \mathbf G(T)$ defines
a modification from the category of latin bitrades to the category
of abelian groups. This also follows immediately from 
\cite[Proposition~3.1]{dk1} since $\mathbf G(T)$ can be identified
with $G(T(\ll))/(G(T(\ll))'$ (by $G$ one denotes in \cite{dk1}
the modification into the category of all groups).

\begin{lem}\label{46}
Let $T=(T^\ll,T^\rr)$ be a spherical latin bitrade.
Then $\mathbf g_T$ provides an embedding of $T(\ll)$ into 
$\mathbf G(T)$.
\end{lem}
\dem Let $a=(a_1,a_2,a_3) \in T^\ll$, $b = (b_1,b_2,b_3) \in T^\ll$
and $i \in \{1,2,3\}$ be such that $a_i \ne b_i$. Put
$\mathbf g_T = (g_1,g_2,g_3)$. Our goal is to prove that 
$g_i(a_i) \ne g_i(b_i)$. By \tref{T5} there exists an abelian
group $K$ and a homotopy $(\sigma_1,\sigma_2,\sigma_3)$ of $T(\ll)$
into $K$ such that $\sigma_i(a_i) \ne \sigma_i(b_i)$. Because
$\mathbf g_T$ is a modification, there must exist a group homomorphism
$\varphi:\mathbf G(T) \to K$ such that $\sigma_j = \varphi g_j$
for all $j \in \{1,2,3\}$. But that means $g_i(a_i) \ne g_i(b_i)$.
\edem

Let again $T=(T^\ll,T^\rr)$ be a (general) latin bitrade.
Put $\mathbf g_T = (g_1,g_2,g_3)$. Following
\cite{dk1} define $\mathbf H(T)$ as the subgroup of $\mathbf G(T)$
generated by the set of all $g_i(b_i) - g_i(a_i)$, where
$(a_1,a_2,a_3),(b_1,b_2,b_3) \in T^\ll$ and $i \in \{1,2,3\}$.
Note that to generate $\mathbf H(T)$ we may consider only $i \in \{1,2\}$, 
and that we also may keep $a=(a_1,a_2,a_3)$ fixed if $b = (b_1,b_2,b_3)$
runs through $T^\ll$ (to see the latter note that 
$g_i(b_i) - g_i(c_i) = (g_i(b_i) - g_i(a_i)) - (g_i(c_i) - g_i(a_i))$).
Put $\mathbf g_T[a] = (h_1,h_2,h_3)$. From \lref{45} we immediately
obtain the following observation:
\begin{lem}\label{47} 
The triple $\mathbf g_T[a]=(h_1,h_2,h_3)$ is a homotopy from $T(\ll)$ to 
$\mathbf H(T)$, and $\mathbf g_T[a]$ embeds $T(\ll)$
into $\mathbf H(T)$ if and only $\mathbf g_T$ embeds 
$T(\ll)$ into $\mathbf G(T)$. Furthermore, $\mathbf H(T)$ is generated
by the set  $\{h_1(b_1),h_2(b_2);$ $(b_1,b_2,b_3) \in T^\ll\}$.
\end{lem}

\begin{thm}\label{48}
Let $T = (T^\ll, T^\rr)$ be a spherical latin bitrade.
Then $T(\ll)$ can be embedded into the abelian group $\mathbf H(T)$
and this group is finite.
\end{thm}
\dem The fact that $T(\ll)$ embeds into $\mathbf H(T)$ follows
immediately from Lemmas~\ref{46} and~\ref{47}. The group $\mathbf
H(T)$ has only finitely many generators and so it suffices to show
that it is a torsion group. Choose $a \in T^\ll$ and assume the contrary.
Then there exists a surjective group homomorphism $\pi: \mathbf H(T) \to
\mathbb Z$, and from \lref{47} we see that $(\pi h_1,\pi h_2, \pi h_3)$ 
is a homotopy $T(\ll) \to
\mathbb Z$ such that $\mathbb Z$ is generated by the set of all
$\pi h_1(b_1)$ and $\pi h_2(b_2)$, where $b = (b_1,b_2,b_3)$ runs
through $T^\ll$. Since $\pi h_i(a_i) = 0$ for all $i \in \{1,2,3\}$
we also see that the homotopy 
$(\pi h_1,\pi h_2, \pi h_3)$ is not trivial. However, that contradicts
\lref{42}.
\edem

Let us mention that not all spherical latin bitrades embed into
a cyclic group---the least known example has size 12 and is mentioned 
already in \cite{drath}. Ian Wanless found in 2006 several latin bitrades
that can be embedded into no group and are of size 24 and genus 4.
The following bitrade is toroidal and also embeds into no group.
It has 5 rows, 6 columns, 7 symbols, and is of size 18.\\

\smallskip
\begin{center}
$
\begin{array}{|c||c|c|c|c|c|c|}
\hline 
\ll & f & a & b & c & d & g \\
\hline
\hline e & 1 &   & 3 & 4 & 5 &  \\
\hline x & 3 & 6 & 2 & 5 &   & 7 \\
\hline y & 5 &   &   &   & 1 &    \\
\hline z & &   & 4 & 2 &   &   \\
\hline t & 7 & 5 & 6 & 3 &   & 2 \\
\hline
\end{array}
$
\hskip 4em
$
\begin{array}{|c||c|c|c|c|c|c|c|}
\hline
\rr & f & a & b & c & d & g \\
\hline
\hline e & 3 &   & 4 & 5 & 1 & \\
\hline x & 7 & 5 & 6 & 3 &   & 2 \\
\hline y & 1 &   &   &   & 5 &   \\
\hline z &  &   & 2 & 4 &   &   \\
\hline t & 5 & 6 & 3 & 2 &   & 7 \\
\hline
\end{array}
$
\end{center}
\vspace{10pt}

Denote the left latin trade as $T(\ll)$ and the right trade $T(\rr)$.
It is not difficult to see that $\mathbf H(T(\rr)) \cong \mathbb Z_{10}$
and that $T(\rr)$ embeds into $\mathbf H(T(\rr))$. In fact 
$H(T(\rr))$, where $H$ means (as in \cite{dk1}) the noncommutative version of 
$\mathbf H$, is isomorphic to $\mathbb Z_{10}$ as well.
Our claim that the bitrade embeds into no group refers to $T(\ll)$.
To see that it embeds into no abelian group
is easy, but the general case seems to deserve a formal proof:
\begin{lem} All homotopies of $T(\ll)$ into a group are trivial.
\end{lem}
\dem The rows are denoted by $e$, $x$, $y$, $z$ and $t$, and the columns
by $f$, $a$, $b$, $c$, $d$ and $g$. We shall regard
these elements as generators of a group $G$ to which there exists
a homotopy from $T(\ll)$. If the homotopy is nontrivial,
then there exists a nontrivial homotopy in which $e = 1 =f$. To
see that modify \lref{46} so that it is valid for noncommutative
groups as well (\cite[Lemma~1.6]{dk2} or \cite[Lemma~3.3]{dk1}). Assuming
$e =f =1$ we get $b=eb = xf = x$, $c = ec = zb = zx$, $d = 
ed = yf = y = xc = xzx$ and $zc = xb = x^2$. Now, $c = zx$
yields $z^2x = x^2$, and we obtain $x = z^2$. Using $t = xc\m$
and $g = x\m t$ we find that all generators are powers of $z$.
In particular, $e = 1$, $x = z^2$, $y = z^5$, $t = z\m$, $f = 1$,
$b = z^2$, $c = z^3$, $d = z^5$, $g = z^{-3}$ and $a = t\m y =
z^6$. Now, $z = tb = xa = z^8$ yields $z^7 = 1$, and so $1 =
ef = y^2 = z^{10}$ implies $z =1$. Group $G$ is thus trivial,
and so there exists no nontrivial homotopy $T(\ll) \to G$.
\edem

The question whether it is possible to embed every spherical
latin bitrade into an abelian group got certain publicity at the
workshop ``Algebraic and geometric aspects of latin trades'' that was
organized at Charles University, Prague, in February 2006.
After submitting the first version of this paper we contacted
Cavenagh and Wanless since we knew that they had been working
upon the problem. It turned out that they found
(amongst others) another proof \cite{cw}. 
Both research efforts have been independent.
There are several
common features, but there are also quite a few dissimilarities---e.g. 
\cite{cw} does not use dissections. 

We finish by two problems. Let $T = (T^\ll, T^\rr)$ be a spherical latin bitrade.

\noindent (1) Is it always possible to retrieve $T$ from a dissection
if $\mathbf H(T)$ is cyclic?

\noindent (2) Must $\mathbf H(T)$ be cyclic when $T$ can be derived
from a separated dissection?

\subsection*{Acknowledgment} We thank to both referees for suggestions
that improved presentation of this paper. They also helped us to
correct notational inconsistencies and small mistakes that were
contained in the submitted manuscript.

\end{document}